\documentclass{amsart}
\usepackage{amssymb}
\usepackage{euscript}
\usepackage{amsfonts}
\usepackage{epsfig}
\usepackage{color}
\usepackage[frame,ps,dvips,matrix,arrow,curve,rotate]{xy}
\let\oldcite\cite

\newtheorem{thm}{Theorem}[section]
\newtheorem{cor}[thm]{Corollary}
\newtheorem{lem}[thm]{Lemma}

\newtheorem{prop}[thm]{Proposition}
\theoremstyle{definition}
\newtheorem{defn}[thm]{Definition}
\theoremstyle{remark}
\newtheorem{rem}[thm]{Remark}
\numberwithin{equation}{section} \theoremstyle{remark}
\newtheorem{ex}[thm]{Example}

\newcommand{\mfG}{\mathfrak{G}}
\newcommand{\mfH}{\mathfrak{H}}
\newcommand{\mfF}{\mathfrak{F}}
\newcommand{\mfK}{\mathfrak{K}}
\newcommand{\mfC}{\mathfrak{C}}
\newcommand{\mfD}{\mathfrak{D}}

\newcommand{\tcat}{\mathbf{2Cat}}

\newcommand{\tgpd}{\mathbf{2Gpd}}
\newcommand{\tgpdw}{\mathbf{2Gpd}_w}
\newcommand{\wtgpd}{\mathbf{W2Gpd}}
\newcommand{\sset}{\mathbf{SSet}}

\let\:=\colon

\newcommand{\lra}{\longrightarrow}
\newcommand{\llra}[1]{\stackrel{#1}{\lra}}

\renewcommand{\hom}{\operatorname{Hom}}
\newcommand{\Hom}{\operatorname{Hom}}

\newcommand{\Hombf}{\mathbf{Hom}}
\newcommand{\homw}{\underline{\mathcal{HOM}}}
\newcommand{\homt}{\underline{\hom}}
\newcommand{\homs}{\underline{\mathcal{H}om}}
\newcommand{\Mapbf}{\underline{\mathcal{RH}om}}

\newcommand{\id}{\operatorname{id}}
\newcommand{\Ho}{\operatorname{Ho}}

\newcommand{\Cosk}{\operatorname{Cosk}}

\def\smashedlongrightarrow{\setbox0=\hbox{$\longrightarrow$}\ht0=1pt\box0}
\def\risom{\buildrel\sim\over{\smashedlongrightarrow}}
\def\smashedst{\setbox0=\hbox{$\rightrightarrows$}\ht0=4pt\box0}
\def\smashedar{\setbox0=\hbox{$\longrightarrow$}\ht0=4pt\box0}
\newcommand{\sst}[1]{\stackrel{#1}{\smashedst}}
\newcommand{\arr}[1]{\stackrel{#1}{\smashedar}}
\newcommand{\oux}[2]{\underset{#1}{\overset{#2}\times}}

\newcommand{\uc}[2]{{#1}_{#2}^{\urcorner}}
\newcommand{\lc}[2]{{}_{\llcorner}\!{#1}_{#2}}

\begin{document}

\title{Notes on 2-groupoids, 2-groups and crossed modules}%

\author{Behrang Noohi}%

\begin{abstract}
    This paper contains some basic results on 2-groupoids, with
    special emphasis on computing derived mapping 2-groupoids between
    2-groupoids and proving their invariance under  strictification.
    Some of the results proven here are presumably
    folklore (but do not appear in the literature to the author's knowledge)
    and some of the results seem to be new. The main technical tool
    used throughout the paper is the Quillen model structure on the
    category of 2-groupoids introduced by Moerdijk and Svensson.
\end{abstract}

\maketitle

\section{Introduction}{\label{S:Intro}}

Crossed modules, invented by J.H.C.\ Whitehead in the early
twentieth century, are algebraic models for (connected) homotopy
2-types (i.e.\ connected spaces with no homotopy group in degrees
above 2), in much the same way that groups are algebraic models for
homotopy 1-types. It is a standard fact that
 a crossed module is essentially the same thing as a 2-group
 (that is, a 2-category with one object in which all 1-morphisms
 and all 2-morphisms are strictly invertible). This is made precise by establishing
 an equivalence between the category of crossed modules and  the category
 of 2-groups (Section~\ref{S:Review}). In particular,   2-groups
 can also be used to model
 connected 2-types. More generally, 2-group{\em oid}s can be used
 to model arbitrary (i.e. not necessarily connected) 2-types.

\enlargethispage{10pt}

But what do we really mean when we say  2-groupoids model homotopy
2-types?

According to Whitehead, to a connected homotopy 2-type one can
associate a crossed module and from this crossed module one can
recover the original 2-type. More generally, to any 2-type one can
associate a 2-groupoid (the Whitehead 2-groupoid) and from this
2-groupoid one can recover the original 2-type. One should keep in
mind, however, that understanding 2-types is not only about
understanding a single space with vanishing $\pi_i$, $i\geq 3$, but
also about understanding {\em mapping spaces} between such spaces.
So one may ask: is the category of 2-groupoids rich enough to
provide us with such information? For instance, given 2-groupoids
$\mfG$ and $\mfH$, can we capture the set of homotopy classes of
maps between the corresponding 2-types from the 2-groupoids of
morphisms between $\mfG$ and $\mfH$?

The answer to the latter question is indeed no. There are usually
far too few 2-groupoid morphisms to produce all homotopy classes. As
a consequence, to have correct algebraic models for mapping spaces
between 2-types, one needs to resort to more sophisticated tools.

There are two standard ways to go about this problem. One is to
enlarge the set of morphisms by taking into account   {\em weak}
morphisms between 2-groupoids. The other is to endow the category of
2-groupoids with a model structure and use a cofibrant replacement
procedure to define {\em derived} mapping spaces. There seems to be
a common understanding among experts how these things should be done
and agreement that both approaches yield essentially the same
results. But to my knowledge, there are not many, if any, places in
the literature where one can look for precise statements, let alone
proofs.

This article is an attempt to fill this gap. The goal is to present
various models for the mapping space between two given 2-groupoids
(or weak 2-groupoids) and to verify that such models are naturally
homotopy equivalent and that they provide algebraic models for the
mapping space between the corresponding 2-types. We will also do the
pointed version of these results. The results of this paper are used
in \cite{Maps} to provide yet another model for the space of weak
morphisms between two crossed modules, which is very simple and
handy and has interesting applications to the study of group actions
on stacks.

This paper is organized as follows. There are essentially three
parts, each one introducing a certain degree of ``weakness'' to
the previous one!

In the first part, which includes up to Section~\ref{S:2}, we are
concerned with {\em strict} 2-groupoids and {\em strict} functors
between them, and we try to understand the derived mapping spaces
(read, 2-groupoids)
 between 2-groupoids. (We will, however, already need to use weak
transformations in this section, as strict transformation are not
homotopically correct  in this context.) Our main technical tool is
the Moerdijk-Svensson model structure, and its monoidal
strengthening developed in \cite{Lack}.

In the second part of the paper, Section~\ref{S:2} and
Section~\ref{S:Crossed2}, we allow weak 2-functors into the picture,
and study the corresponding simplicial mapping spaces. The point is
that allowing this extra flexibility lets us get away with not
having  to make cofibrant replacements when computing derived
mapping spaces. That is, the 2-groupoid of weak 2-functors, weak
2-transformations and modification between arbitrary 2-groupoids
$\mfH$ and $\mfG$ has the correct homotopy type and is homotopy
equivalent to the derived mapping space.

In the third part of the paper,  Section~\ref{S:3}, we bring in weak
2-groupoids. We see that  introducing weak 2-groupoids essentially
does not add anything to the homotopy theory. Put differently, we
can always replace weak 2-groupoids with strict ones without losing
any homotopy theoretic information, and this strictification does
not alter derived mapping spaces.

Since the author's ultimate goal is to apply these results to
2-groups and crossed modules, we  also develop a pointed version of
the theory, parallel to the unpointed one, and explain briefly how
the pointed theory can be translated to the language of crossed
modules.

\vspace{0.1in} \noindent{\em Acknowledgement.}  I would like  to
thank the Max-Planck-Institut f\"{u}r Mathematik, Bonn, where the
research for this work was done, for providing pleasant working
conditions. I thank Joachim Kock for reading the manuscript and
making many useful remarks, Fernando Muro for many useful
discussions and for explaining Andy Tonks' results to me, and
Bertrand To\"en for pointing out the work of Steve Lack~\cite{Lack}.
I would also like to thank Tim Porter for useful comments on an
earlier  version of this paper. Thanks are due to an anonymous
referee for suggesting inclusion of the discussion of monoidal
properties of the functors $N$ and $W$ in $\S$~\ref{SS:LaxMon}.

\tableofcontents

\section{Notations and terminology}{\label{S:NT}}

A {\em 2-functor} between 2-categories means a strict 2-functor.
 We sometimes refer to a 2-functor simply as a functor, or, more informally,
  a {\em map of 2-categories}. A
 2-functor between 2-groupoids is  simply a 2-functor between the
 underlying 2-categories.

In this paper, we encounter various kinds of ``hom''s.  To make it
easier to keep track of these we follow certain conventions. First,
the underlined ones are internal, namely they are 2-groupoids or
2-categories. The more curly they get, the more weakness is
involved. E.g., $\underline\Hom$ means the 2-category of strict
2-functors, strict 2-transformations and modifications. $\homw$ then
stands for the 2-category of weak 2-functors, weak transformations
and modifications. $\Hombf$ is used for simplicial mapping space in
the category of simplicial sets.

 For objects $A$ and $B$ in a  category $\mathbf{C}$ with a notion
 of weak equivalence, we denote the set  of morphisms
 in the homotopy category from $B$
 to $A$, that is
$\Hom_{\Ho(\mathbf{C})}(B,A)$, by $[B,A]_{\mathbf{C}}$.

 The terms `morphism', `1-morphism' and `arrow' will be used synonymously.
 We use multiplicative notation for 1-morphisms
 of a 2-category, as opposed to the compositional notation (that is, $fg$
 instead of $g\circ f$). We also do the same with the horizontal
 and vertical
 composition of 2-morphisms. To avoid confusion, we put our
 2-morphisms inside square brackets when we are doing vertical
 compositions. This is all illustrated in the following example:
 
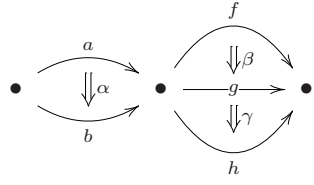
\[\xymatrix@C=13pt@R=-2pt@M=6pt{&  & &
   \ar@{=>}[dd]^{\beta}&\\  &  \ar@{=>}[dd]^{\alpha} &&&\\
       \bullet \ar@/^1pc/ [rr]^{a} \ar@/_1pc/ [rr]_{b} &&
       \bullet  \ar@/^2pc/ [rr]^{f} \ar[rr] |-{g} \ar@/_2pc/ [rr]_{h} &
                    \ar@{=>}[dd]^{\gamma} &  \bullet    \\ &&&& \\
                                                          &&   & &   }\]

For example, the 2-morphism from $af$ to $ag$ is denoted by
$a\beta$, and the 2-morphism from $af$ to $bh$ is given by
 \[[a\beta][a\gamma][\alpha h]=[a\beta][\alpha h][ b\gamma]=
  [\alpha f][b\beta][b\gamma]=[\alpha\beta][b\gamma]=[a\beta][\alpha\gamma].\]

 We use the German letters $\mfC$, $\mfD,\ldots$ for general 2-categories
 and $\mfG$, $\mfH,\ldots$ for 2-groupoids. The
 upper case Latin letters $A$, $B$, $C,\ldots$ are used
 for objects in such 2-categories,
 lower case script letters $a$, $b$, $g$, $h,\ldots$ for 1-morphisms, and
 lower case Greek  letters $\alpha$, $\beta,\ldots$ for 2-morphisms.

The symbol $\simeq$ is used for weak or homotopy equivalences, and
$\cong$ is used for isomorphisms.

We denote by $\tgpd$ the category of small 2-groupoids and strict
functors, by $\tgpdw$ the category  of small 2-groupoids and weak
functors (Section~\ref{S:2}), and by $\wtgpd$ the category of small
weak
 2-groupoids and weak 2-functors (Section~\ref{S:3}). {\em Our
notions of weak 2-functor and weak 2-category are slightly stronger
than the standard notions in that we do not weaken the identities.}
The category of simplicial sets is denoted by $\sset$, and the
category of small 2-categories by $\tcat$. Various nerve functors
that appear in this paper are all denoted by $N$. The Whitehead
functor $W \: \sset \to \tgpd$, to be introduced in
Section~\ref{S:M-S}, is the left adjoint to $N \: \tgpd \to \sset$.

 By {\em fiber product} of 2-ca\-te\-go\-ries we mean strict fi\-ber pro\-duct.
 We will not en\-coun\-ter the {\em homotopy} fiber pro\-duct of 2-ca\-te\-go\-ries
 in this pa\-per.

 For notation concerning crossed modules see
 $\S$\ref{SS:CrossedMod}.

\section{Quick review of 2-groups and crossed modules}{\label{S:Review}}

\subsection{Quick review of 2-groups}{\label{SS:2gp}}

We recall some basic facts about 2-groups and crossed modules. Our
main references are \cite{Baez,Baues,Brown,Loday,M-S,M-W,W}. The
reader is also encouraged to consult the works of the Bangor school
(mostly available on their webpage~\cite{Bangor}).

 A {\bf 2-group} $\mathfrak{G}$ is a  group object in the category of groupoids.
 Alternatively, we can define a 2-group to be a groupoid object in the category
 of groups, or  also, as a (strict) 2-category with one object in
 which all  1-morphisms and 2-morphisms are invertible (in the strict sense).
 We will try to adhere to the `group in groupoids' point of view throughout the paper, but
 occasionally switching back and forth between different points of view
 is inevitable. Therefore, the reader will  find it rewarding to master how
 the equivalence of these three point of views works.

 A (strict) {\em morphism}
 $f \: \mathfrak{G} \to \mathfrak{H}$ of 2-groups is a map of groupoids that
 respects the group operation. If we view $\mathfrak{G}$ and  $\mathfrak{H}$
 as 2-categories with one object, such $f$ is   nothing but  a strict 2-functor.
 The category of 2-groups is denoted by $\mathbf{2Gp}$.

 To a 2-group $\mathfrak{G}$ we associate the groups
 $\pi_1\mathfrak{G}$ and $\pi_2\mathfrak{G}$ as follows. The group
 $\pi_1\mathfrak{G}$
 is the set of isomorphism classes of objects of the groupoid $\mathfrak{G}$.
 The group structure on $\pi_1\mathfrak{G}$ is induced from the group
 operation of $\mathfrak{G}$. The group
 $\pi_2\mathfrak{G}$ is the group of automorphisms  of the identity
 object $e \in \mathfrak{G}$. This is an abelian group.
 The groups $\pi_1$ and $\pi_2$ are functorial in $\mathfrak{G}$.
 A  morphism
 $f \: \mathfrak{G} \to \mathfrak{H}$ of 2-groups is called an {\em equivalence}
 if it induces  isomorphisms on $\pi_1$ and $\pi_2$. The {\em homotopy
 category} of 2-groups is the category obtained by inverting all the
 equivalences in $\mathbf{2Gp}$. We denote it by
 $\Ho(\mathbf{2Gp})$.

 \vspace{0.1in}
 \noindent {\em Caveat}: an equivalence between 2-groups
 need not have an inverse. Also, two equivalent
  2-groups may not be related by an equivalence, but only a zig-zag of
  equivalences.

 \subsection{Quick review of crossed modules}{\label{SS:CrossedMod}}

 A {\bf crossed module}
 $\mathfrak{G}=[\varphi \: G_2 \to G_1]$
  is a pair of groups $G_1,G_2$, a group homomorphism $\varphi \: G_2 \to G_1$,
  and a (right) action of $G_1$ on $G_2$, denoted $-^a$,
  which lifts the conjugation action
  of $G_1$ on
  the image of $\varphi$ and descends the conjugation action of $G_2$ on itself.
  Before making this precise, let us introduce some (perhaps
  non-standard) notation.

 \vspace{0.1in}
 \noindent{\bf Notation.}   We usually denote the elements of $G_2$ by
Greek letters and those of $G_1$ by lower case Roman letters.  We
sometimes suppress $\varphi$ from the notation
 and denote $\varphi(\alpha)$ by $\underline{\alpha}$.  For elements
 $g$ and $a$ in a group $G$ we sometimes denote  $a^{-1}ga$
 by $g^a$.

\vspace{0.1in}

 With the above notation,
 the two compatibility axioms of a crossed module can be written in the
 following way:
 \begin{itemize}
   \item[$\mathbf{CM1}.$] \ $\forall \alpha,\beta \in G_2, \
                                \beta^{\underline{\alpha}}=\beta^{\alpha}$;
   \item[$\mathbf{CM2}.$] \ $\forall \beta \in G_2, \forall a\in G_1, \
                                   \underline{\beta}^a=\underline{\beta^a}$.
 \end{itemize}

  The kernel of $\varphi$ is a central (in particular abelian) subgroup of $G_2$
  and is denoted by $\pi_2\mathfrak{G}$. The image of $\varphi$ is a normal
  subgroup of $G_1$ whose cokernel is denoted by  $\pi_1\mathfrak{G}$.
  A  (strict) morphism  $P: \mfH \to \mfG$ of crossed modules is a pair of
  group homomorphisms
  $p_2 \: H_2 \to G_2$ and $p_1 \: H_1 \to G_1$ which commute
  with the $\varphi$ maps and respect the actions. Such a
  morphism induces homomorphisms on $\pi_1$ and $\pi_2$.
  A  morphism is called an {\em equivalence} if
  it induces isomorphisms on $\pi_1$ and $\pi_2$.

 \subsection{Equivalence of 2-groups and crossed modules}{\label{SS:Equiv}}
   There is a natural pair of inverse equivalences   between
   the category $\mathbf{2Gp}$ of 2-groups  and the category
   $\mathbf{CrossedMod}$ of crossed modules.
   Furthermore,
   these functors preserve $\pi_1$ and $\pi_2$. They are constructed
   as follows. We only describe the effect of these functors on
   objects.

   \medskip
   \noindent{\em Functor from 2-groups to crossed modules.}
   Let $\mathfrak{G}$ be a 2-group. Let $G_1$ be the group of   objects of
   $\mathfrak{G}$, and $G_2$  the set of   arrows emanating from the
   identity object $e$; the latter is also a group (namely, it is a subgroup
   of the group  of arrows of $\mathfrak{G}$).

   Define the map $\varphi \: G_2 \to G_1$ by sending $\alpha \in G_2$ to $t(\alpha)$.

   The action of $G_1$ on $G_2$ is given by conjugation. That is,
   given $\alpha \in G_2$ and $g \in G_1$, the action is given by
   $g^{-1}\alpha g$. Here were are thinking of $g$ as an identity
   arrow and multiplication takes place in the group of arrows of
   of $\mfG$. It is readily checked that $[\varphi \: G_2 \to G_1]$
   is a crossed module.

    \medskip
   \noindent{\em Functor from crossed modules to 2-groups.}
   Let  $[\varphi \: G_2 \to G_1]$ be a crossed module. Consider the
   groupoid $\mathfrak{G}$ whose underlying set of objects is $G_1$ and whose
   set of arrows is $G_1\rtimes G_2$. The source and target maps are given by
   $s(g,\alpha)=g$, $t(g,\alpha)=g\varphi(\alpha)$.

   Recall that the group operation in $G_1\rtimes G_2$ is defined
   by $(g,\alpha)(h,\beta)=(gh,\alpha^h\beta)$.  It is easy to
   see that $s$, $t$, and the inclusion map $G_1 \hookrightarrow G_1\rtimes
   G_2$ are group homomorphisms. So $\mfG$ is naturally a group
   object in the category of groupoids.

 \medskip The above discussion shows that there is a pair of
inverse functors inducing an equivalence between
$\mathbf{CrossedMod}$ and $\mathbf{2Gp}$. These functors respect
$\pi_1$ and $\pi_2$. Therefore, we have an equivalence
   \[\xymatrix@=16pt@M=8pt{
   \Ho(\mathbf{CrossedMod})  \ar@<0.5ex>[r]
                 & \ar@<0.5ex>[l] \Ho(\mathbf{2Gp}).}\]
This equivalence respects the functors $\pi_1$ and $\pi_2$.

\section{2-categories and 2-groupoids}{\label{S:1}}

 In this section we quickly go over some basic facts and constructions
 for 2-categories, and fix some terminology.

 By a {\em 2-category} we mean a strict 2-category. A {\em  2-groupoid}
 is a 2-category in which every 1-morphism and every 2-morphism
 has an inverse (in the strict sense). Every category (respectively, groupoid)
 can be thought of as a 2-category (respectively, 2-groupoid) in which
 all 2-morphisms are  identity.

 A {\em 2-functor} between 2-categories means a strict 2-functor.
 We sometimes refer to a 2-functor simply by a functor, or a {\em
 map of 2-categories}. A
 2-functor between 2-groupoids is  simply a 2-functor between the
 underlying 2-categories.

\subsection{Strict transformations between 2-functors}{\label{SS:Transformations}}

Let $\mathfrak{C}$ and  $\mathfrak{D}$ be (strict) 2-categories, and
let $P,Q \: \mathfrak{D} \to \mathfrak{C}$ be (strict) 2-functors.
By a {\em 2-transformation} $T \: P \Rightarrow Q$ we mean an
assignment  of an arrow $t_A \: P(A) \to Q(A)$ in $\mathfrak{C}$ to
every object $A$ in $\mathfrak{D}$ such that for every 2-morphism
$\gamma\: c \Rightarrow c'$ in $\mathfrak{D}$ between two arrows $c$
and $c'$ with common source and targets $A$ and $B$,  we have
$P(\gamma)t_B=t_AQ(\gamma)$. This can be depicted as the following
commutative pillow:
    \[\xymatrix@=8pt@M=10pt{
         P(A) \ar[rr]^{t_A} \ar@/_1.1pc/ [dd]_{P(c)}  \ar@{}[d]
                            |(0.8){P(\gamma)} \ar@/^1.1pc/ [dd]^{P(c')} &&
              Q(A)  \ar@/_1.1pc/ [dd]_{Q(c)} \ar@{}[d]
                            |(0.8){Q(\gamma)}  \ar@/^1.1pc/ [dd]^{Q(c')}   \\
           \Longrightarrow &&  \Longrightarrow \\
         P(B)    \ar[rr]_{t_{B}}        &&  Q(B)   }\]  \label{pillow}

Alternatively, a 2-transformation $P \Rightarrow Q$ is the same as a
2-functor $\mathbf{I}\times\mfD \to \mfC$ whose restriction to
$\{0\}\times\mfD$ is $P$ and whose restriction to $\{1\}\times\mfD$
is $Q$. Here $\mathbf{I}$ stands for the category $\{0 \to 1\}$. We
sometimes refer to a 2-transformation simply by a transformation, if
there is no chance of confusion.

 A transformation between two 2-transformations $S$, $T$,  called
 a {\em modification}, is a rule assigning to each object $A \in  \mfD$
 a 2-morphism $\mu_A \: s_A \Rightarrow t_A$ in  $\mfC$, such that
 for every arrow $c \: A \to B$ in $\mfD$, we have
 $P(c)\mu_B=\mu_AQ(c)$. This can be depicted as the following commutative pillow:
  \[\xymatrix@C=13pt@R=-4pt@M=6pt{ & \ar@{=>}[dd]^{\mu_A}&\\
        P(A) \ar[dddddd]_{P(c)} \ar@/^1pc/ [rr]^{s_A} \ar@/_1pc/ [rr]_{t_A} &&  Q(A)  \ar[dddddd]^{Q(c)} \\
            & &  \\ \\ \\ \\
    & \ar@{=>}[dd]^{\mu_B}&\\
        P(B) \ar@/^1.1pc/ [rr]^{s_B} \ar@/_1.1pc/ [rr]_{t_B} &&  Q(B)   \\
            & &   }\]

Alternatively, a modification $S \Rightarrow T$ is the same as a
2-functor $\mathbb{I}\times\mfD \to \mfC$ whose restriction to the
subcategory $(s,\id)\: \mathbf{I}\times \mfD \hookrightarrow
\mathbb{I}\times \mfD$ corresponds to the transformation $S$, and
whose restriction to the subcategory $(t,\id)\: \mathbf{I}\times
\mfD \hookrightarrow \mathbb{I}\times \mfD$ corresponds to the
transformation $T$. Here $\mathbb{I}$ stands for the 2-category
 \[\xymatrix@C=10pt@R=-4pt@M=6pt{ & \ar@{=>}[dd]^{\mu}&\\
        0 \ar@/^1.1pc/ [rr]^{s} \ar@/_1.1pc/ [rr]_{t} &&  1   \\
            & &    }\]

\begin{defn}{\label{D:homo1}}
  Given 2-categories $\mathfrak{C}$ and $\mathfrak{D}$,
  we define $\homt(\mathfrak{D},\mathfrak{C})$ to be the
  2-category whose objects are 2-functors from
  $\mathfrak{D}$ to $\mathfrak{C}$, whose 1-morphisms
  are 2-transformations between 2-functors, and whose
  2-morphisms are modifications.
\end{defn}

 The following exponential law is easy to prove.

\begin{lem}{\label{L:exponential}}
   Let $\mathfrak{C}$, $\mathfrak{D}$ and $\mathfrak{E}$
   be 2-categories, then we have a natural isomorphism
   of 2-categories
    \[\homt(\mathfrak{E} \times \mathfrak{D}, \mathfrak{C})
        \cong
         \homt(\mathfrak{E},\homt(\mathfrak{D},\mathfrak{C})).\]
\end{lem}

 With this, $\mathbf{2Cat}$ is enriched in
2-categories. When $\mathfrak{C}$ is a 2-groupoid,
$\homt(\mathfrak{D},\mathfrak{C})$ is also a 2-groupoid. So,
$\mathbf{2Gpd}$ is naturally enriched in 2-groupoids. In other
words, $\mathbf{2Cat}$ and $\mathbf{2Gpd}$ are closed monoidal
categories, the product being cartesian product.

\subsection{Weak transformations between 2-functors}{\label{SS:Transformations2}}

There is also a notion of {\em weak} 2-transformation between
(strict) 2-functors, and modifications between them. We recall the
definitions.

 A {\em weak
transformation} between  two 2-functors $P$ and $Q$ is
 like a strict transformation, except that to a
  morphism  $c \: A\to B$ in $\mfD$ we assign
 a 2-cell
 \[\xymatrix@=12pt@M=10pt{
         P(A) \ar[r]^{t_A}\ar[d]_{P(c)}  &
              Q(A)  \ar[d]^{Q(c)}   \\
         P(B)    \ar[r]_{t_{B}}   \ar@{=>}[ur]^{\theta_c}     &  Q(B)   }\]
We require that $\theta_{\id}=\id$, and that $\theta_h$
 satisfy the obvious compatibility conditions with respect to
2-morphisms and composition of morphisms. There are two types of
conditions here. One is the same as the commutative pillow of page
\pageref{pillow}, in which the top and bottom squares  are now
decorated with $\theta_c$ and $\theta_{c'}$, respectively. (Just to
make sure the reader visualizes the pillow correctly, the top square
of the pillow is exactly the 2-commutative square of the previous
paragraph.) The second compatibility condition says that
 for every two
composable morphisms $a$ and $b$, the corresponding squares
compose. By abuse of notation, this can be written as
$\theta_a\theta_b=\theta_{ab}$.

Modifications between two weak 2-transformations are defined in
the same way as they were defined for strict 2-transformations
($\S$\ref{SS:Transformations}).

\begin{defn}{\label{D:homo3}}
  Given 2-categories $\mathfrak{C}$ and $\mathfrak{D}$,
  we define $\homs(\mathfrak{D},\mathfrak{C})$ to be the
  2-category whose objects are 2-functors from
  $\mathfrak{D}$ to $\mathfrak{C}$, whose 1-morphisms
  are weak 2-transformations between 2-functors, and whose
  2-morphisms are modifications.\footnote{In \oldcite{Gray},
  this 2-category is denoted by $\operatorname{Fun}$.}
\end{defn}

Of course, the exponential law of Lemma~\ref{L:exponential} will not
work with $\homs$ unless we change our product accordingly. It turns
out that this can be done, namely by replacing the usual product
with the Gray tensor product. For a nice description of Gray tensor
product see \cite{Lack2}, last paragraph of Section 2.

In the following, $\otimes$ stands for the Gray tensor product.

\begin{lem}[\oldcite{Gray}, page 73, Theorem I,4.9]{\label{L:Grayexponential}}
   Let $\mathfrak{C}$, $\mathfrak{D}$ and $\mathfrak{E}$
   be 2-categories, then we have a natural isomorphism
   of 2-categories
    \[\homs(\mathfrak{E} \otimes \mathfrak{D}, \mathfrak{C})
        \cong
         \homs(\mathfrak{E},\homs(\mathfrak{D},\mathfrak{C})).\]
\end{lem}

This gives us another enrichment of $\mathbf{2Cat}$ over itself.
When $\mathfrak{C}$ is a 2-groupoid, then
$\homs(\mathfrak{D},\mathfrak{C})$ is also a 2-groupoid. So,
similarly, we have a new enrichment of $\mathbf{2Gpd}$ over
itself. Therefore, each of $\tgpd$ and $\tcat$ admits two
different closed monoidal structures, one with $\homt$ and
cartesian product, the other with $\homs$ and Gray tensor product.

\subsection{Nerve of a 2-category}{\label{SS:Nerve}}

We review the nerve construction for 2-categories, and recall its
basic properties \cite{M-S}.

Let $\mfG$ be a 2-category .  We define the {\em nerve} of $\mfG$, denoted by
$N\mfG$, to be the simplicial set defined as follows.
 The set  of 0-simplices
of $N\mfG$ is the set of objects of $\mfG$. The 1-simplices are
the morphisms in $\mfG$. The 2-simplices are diagrams of the form
\[\xymatrix@C=6pt@R=14pt@M=6pt{  & B \ar[dr]^g \ar@{=>}[d]^{\alpha} &   \\
            A \ar[rr]_h \ar[ru]^f   &    &   C   }\]
where $\alpha \: fg \Rightarrow h$ is a 2-morphism. The
3-simplices of $N\mfG$ are commutative tetrahedra of the form
\label{nervediagram} \label{3cell}
   \[\xymatrix@C=3pt@R=12pt@M=6pt{ &&& D   &&& \\
                               &&&&  && \\
                            &&& B \ar[uu]^l \ar@{=>}[ul]^{\beta}
                                \ar@{=>}[d]^{\alpha} \ar[drrr]^g &     \ar@{=>}[ul]_{\gamma}  &&   \\
                        A \ar[rrrrrr]_h \ar[urrr]^(0.65)f   \ar[rrruuu]^k
                      && \ar@{:>}[ul]_<{\delta} &&&& C \ar[llluuu]_m   }\]

Commutativity of the above tetrahedron means
$[f\gamma][\beta]=[\alpha m][\delta]$. That is, the following
square of transformations is commutative:
            \[\xymatrix@=12pt@M=10pt{
                 fgm  \ar@{=>}[r]^{f\gamma} \ar@{=>}[d]_{\alpha m}
                                          &  fl \ar@{=>}[d]^{\beta}  \\
                            hm     \ar@{=>}[r]_{\delta}        &    k   }\]
For $n\geq 3$, an $n$-simplex of $N\mfG$ is an $n$-simplex such
that each of its sub 3-simplices is a commutative tetrahedron as
described above. In other words, $N\mfG$ is the coskeleton  of the
3-truncated simplicial set $\{N\mfG_0,N\mfG_1,N\mfG_2,N\mfG_3\}$
defined above.

The nerve  gives us a functor $N \: \mathbf{2Cat} \to
\mathbf{SSet}$, where    $\mathbf{SSet}$ is the the category of
simplicial sets.

\section{Closed model structure on
$\mathbf{2Gpd}$}{\label{S:M-S}}

We quickly review the Moerdijk-Svensson closed model structure on
the category of 2-groupoids. The main reference is \cite{M-S}. A
generalization of this model structure to the case of
crossed-complexes can be found in \cite{BrGo}.

 \begin{defn}{\label{D:2fibration}}
     Let $\mfH$ and $\mfG$ be 2-groupoids, and $P \: \mfH \to \mfG$
     a functor between them. We say that
     $P$ is  a {\em fibration} if it satisfies the
     following properties:
     \begin{itemize}
   \item[$\mathbf{F1.}$] For every arrow $a \: A_0 \to A_1$ in
     $\mfG$, and every object  $B_1$ in $\mfH$ such that $P(B_1)=A_1$,
     there is an object $B_0$ in $\mfH$ and an arrow $b \: B_0 \to B_1$
     such that $P(b)=a$.

   \item[$\mathbf{F2.}$] For every 2-morphism
       $\alpha \: a_0 \Rightarrow a_1$ in
     $\mfG$ and every arrow $b_1$ in $\mfH$ such that $P(b_1)=a_1$,
     there is an arrow $b_0$ in $\mfH$ and a
     2-morphism $\beta \: b_0 \Rightarrow b_1$
     such that $P(\beta)=\alpha$.
  \end{itemize}
 \end{defn}

 \begin{defn}{\label{D:homotopygroup}}
     Let $\mfG$ be a 2-groupoid, and $A$ an object in $\mfG$.
     We define the following.
     \begin{itemize}
   \item $\pi_0\mfG$  is the set of isomorphism classes of objects in $\mfG$.

   \item $\pi_1(\mfG,A)$ is the group of 2-isomorphism classes of arrows
         from $A$ to itself. The {\em fundamental groupoid}
         $\Pi_1\mfG$ is the groupoid obtained by taking the
         groupoid of all 1-morphisms of $\mfG$ and identifying 2-isomorphic
         1-morphisms.

   \item $\pi_2(\mfG,A)$ is the group of 2-automorphisms of the identity
         arrow $1_A \: A \to A$.
  \end{itemize}
  These invariants are functorial with respect to 2-functors.
  A map $\mfH \to \mfG$ is called a {\em (weak) equivalence of 2-groupoids}
  if it induces  bijections on $\pi_0$, $\pi_1$ and $\pi_2$ for every
  choice of basepoint. We usually drop the adjective `weak' when we talk about
  weak equivalences of 2-groupoids, groupoids, 2-groups, or
  crossed modules.
 \end{defn}

We remark that a 2-group in the sense of $\S$\ref{SS:2gp} is
precisely a 2-groupoid with one object. This
  identifies $\mathbf{2Gp}$ as a full subcategory
  of $\mathbf{2Gpd}$. Under this identification,
the notions of $\pi_1$, $\pi_2$, and weak equivalence introduced in
$\S$\ref{SS:2gp} coincide with the ones given above.

Having defined fibrations and weak equivalences between 2-groupoids,
we define {\em cofibrations} using the left lifting property.

\begin{thm}[\oldcite{M-S}, Theorem 1.2]{\label{T:QuillenStr}}
  With weak equivalences, fibrations and cofibrations defined as above,
  the category of 2-groupoids is a closed model
  category.
\end{thm}

In this model structure, every object is fibrant. Cofibrant objects
are a bit trickier, but, morally, they should
be thought of as some sort of {\em free} objects. A more explicit
description of cofibrations in the Moerdijk-Svensson model structure
is given in the Remark on page 194 of \cite{M-S}.

The nerve functor is a bridge between the homotopy theory of
2-groupoids and the homotopy theory of simplicial sets. To justify
this statement, we quote the following from \cite{M-S}.

\begin{prop}[see \oldcite{M-S}, Proposition 2.1]\label{P:nerve}\indent\par
   \begin{itemize}
   \item[$\mathbf{i.}$] The functor $N\: \mathbf{2Cat} \to \mathbf{SSet}$
    sends transformations between 2-functors to simplicial homotopies,
    is faithful and  preserves fiber products. (We will see in
    Section~\ref{S:2} that, in contrast to the case of ordinary categories, $N$
    is not full.)

   \item[$\mathbf{ii.}$] The functor $N \:\mathbf{2Gpd} \to \mathbf{SSet}$
      sends a fibration
      between 2-groupoids (Definition~\ref{D:2fibration}) to a
      Kan fibration. Conversely, if a map $P$ of 2-groupoids is
      sent to a Kan fibration, then $P$ itself is a fibration.
     The nerve of  every
      2-groupoid  is a Kan complex.

   \item[$\mathbf{iii.}$] For every (pointed) 2-groupoid $\mfG$ we have
         $\pi_i(\mfG)\cong\pi_i(N\mfG)$, $i=0,1,2$.

   \item[$\mathbf{iv.}$] A map $f \: \mfH \to \mfG$ of 2-groupoids is an
      equivalence if and only if $Nf \: N\mfH \to N\mfG$ is a weak
      equivalence of simplicial sets.
  \end{itemize}
\end{prop}

 \begin{rem}{\label{R:nerve2gp}}
  Regarding $\mathbf{2Gp}$ as a full subcategory
  of $\mathbf{2Gpd}$, we can also talk about nerves of 2-groups.
  This is a functor $N \: \mathbf{2Gp} \to \mathbf{SSet}_*$,
  where $\mathbf{SSet}_*$ is the category of pointed simplicial sets.
  The above proposition remains valid if we replace 2-groupoids by
  2-groups and $\mathbf{SSet}$ by $\mathbf{SSet}_*$ throughout.
\end{rem}

The  functor $N\: \mathbf{2Gpd} \to \mathbf{SSet}$ has a left
adjoint $W \: \mathbf{SSet} \to \mathbf{2Gpd}$, called the {\em
Whitehead 2-groupoid}, which is defined as follows (see \cite{M-S}
page 190, Example 2). Let $X$ be a simplicial set. The underlying
groupoid of $W(X)$ is $\Pi_1(|X^{(1)}|,|X^{(0)}|)$, where
$|X^{(i)}|$ stands for the geometric realization of the $i^{th}$
skeleton of $X$, and $\Pi_1$ stands for fundamental groupoid. A
2-morphism in $W(X)$ is the equivalence class of a continuous map
$\alpha \: I\times I \to |X|$ such that $\alpha(\{0,1\}\times I)
\subseteq |X^{(0)}|$ and $\alpha(I\times\{0,1\}) \subseteq
|X^{(1)}|$. Two such maps $\alpha$ and $\beta$ are equivalent if
there is homotopy $H\: (I\times I)\times I \to |X|$ between them
such that $H\big((\{0,1\}\times I)\times I\big) \subseteq |X^{(0)}|$
and $H\big((I\times\{0,1\})\times I\big) \subseteq |X^{(1)}|$.

It is easy to see that $W$ preserves homotopy groups. In particular,
it sends weak equivalences of simplicial sets to equivalences
of 2-groupoids. Much less obvious is the following

\begin{thm}[\oldcite{M-S},  Section 2]\label{T:QuillenEq}
   The pair
     \[W\:\mathbf{SSet} \rightleftharpoons \mathbf{2Gpd}:N\]
   is a Quillen pair. It satisfies the following properties:
   \begin{itemize}
   \item[$\mathbf{i.}$] Each adjoint preserves weak equivalences.

   \item[$\mathbf{ii.}$] For every 2-groupoid $\mfG$,  the counit
       of adjunction $WN(\mfG) \to \mfG$ is a weak equivalence.

   \item[$\mathbf{iii.}$] For every simplicial set $X$ such that $\pi_iX=0$,
   $i\geq 3$, the unit of adjunction
      $X \to NW(X)$ is a weak equivalence.
  \end{itemize}
   In particular, the functor
   $N \: \Ho(\mathbf{2Gpd}) \to \Ho(\mathbf{SSet})$ induces an
   equivalence of categories between $\Ho(\mathbf{2Gpd})$ and
   the category of homotopy 2-types. (The latter
   is defined to be the full subcategory of $\Ho(\mathbf{SSet})$
   consisting of all $X$ such that $\pi_iX=0$, $i \geq 3$.)
\end{thm}

\begin{rem}{\label{R:pointed}}
  The pointed version of the above
  theorem is also true.  The proof is just a
  minor modification of the proof of the above theorem.
\end{rem}

It is also well-known that the geometric realization functor
$|-| \: \mathbf{SSet}_* \to \mathbf{Top}_*$     induces an equivalence of
of homotopy categories. So we have  the following

\begin{cor}{\label{C:equiv}}
  The functor $|N(-)| \: \Ho(\mathbf{2Gp}) \to \Ho(\mathbf{Top}_*)$
  induces an equivalence between $\Ho(\mathbf{2Gp})$
   and the full subcategory of $\Ho(\mathbf{Top}_*)$ consisting
   of connected pointed homotopy 2-types.
\end{cor}

\subsection{Monoidal closed model structure on
$\mathbf{2Gpd}$}{\label{SS:Simplicial}}

We saw in Section~\ref{S:M-S} that $\mathbf{2Gpd}$ admits a closed
model structure. We also know from $\S$\ref{SS:Transformations} and
$\S$\ref{SS:Transformations2} that $\mathbf{2Gpd}$ can be made into
a monoidal category in two different ways. Namely, there are two
reasonable ways of enriching the hom-sets to 2-groupoids. In many
applications, it is desirable to have {\em simplicially} enriched
hom-sets. In other words, one would like to have a simplicial
structure in the sense of (\cite{Jardine}, $\S$ II.2, Definition
2.1). Furthermore, to have a full-fledged homotopy theory one also
requires a compatibility between the simplicial structure and the
closed model structure. This compatibility condition is neatly
encoded in Quillen's $\mathbf{SM7}$, (\cite{Jardine}, $\S$ II.3,
Definition 3.1).

We have two candidates for simplicially enriching $\mathbf{2Gpd}$.
We could either take $N\homt(\mfH,\mfG)$, or $N\homs(\mfH,\mfG)$.
It can be shown that $\mathbf{2Gpd}$ does become a simplicial
category with $N\homt(\mfH,\mfG)$. However, $\mathbf{SM7}$ will
not be satisfied, because it is easy to produce examples where the
(derived) mapping spaces do not have the right homotopy type (or
one can directly verify that $\mathbf{SM7}$ fails). At first look,
$N\homs(\mfH,\mfG)$ seems to be the correct alternative, because
in this setting the (derived) mapping spaces can be shown to have
the right homotopy type. However this attempt fails even more
miserably, as this simplicial enrichment is not even well-defined:
given three groupoids $\mfG$, $\mfH$ and $\mfK$, the composition
map
 \[N\homs(\mfK,\mfH)\times N\homs(\mfH,\mfG) \to
  N\homs(\mfK,\mfG)\]
can not even be defined! This is because if we have an arrow in
$\homs(\mfK,\mfH)$, namely a weak 2-transformation between the
2-functors $P_1,Q_1 \:\mfK \to \mfH$, and an arrow in
$\homs(\mfH,\mfG)$, namely a weak 2-transformation between the
2-functors $P_2,Q_2 \:\mfH \to \mfG$, then it is not possible to
define a canonical weak transformation between the 2-functors
$P_2\circ P_1, Q_2\circ Q_1 \: \mfK \to \mfG$ out of these.

In \cite{Lack} it is proven that the monoidal structure on $\tgpd$
defined via $\homs$ and Gray tensor product is compatible with the
Moerdijk-Svensson model structure. Namely, it equips $\tgpd$ with
the structure of a monoidal  closed model category in the sense of
(\cite{Lack}, $\S$ 7). The nerve functor $N \: \tgpd \to \sset$,
however, does not  respect the monoidal structure.

Nevertheless, with some extra work, one can prove that $N$ respects
{\em derived} mapping spaces (Proposition~\ref{P:mappingspace3}).

\begin{defn}{\label{D:derived}}
Let  $\mathbf{C}$ be a monoidal closed model category, with
internal hom-objects denoted by $\homs$. Let
 $A$ and $B$ be objects in
$\mathbf{C}$. The {\em derived  hom-object} $\Mapbf(B,A)$ is
defined to be $\homs(\bar{B},\bar{A})$, where $\bar{A}$ is a
fibrant replacement for $A$ and $\bar{B}$ is a cofibrant
replacement for $B$.
\end{defn}

The derived  hom-object should be thought of as the monoidal
enrichment of  $[B,A]_{\mathbf{C}}$, in the sense that
$\pi_0\Mapbf(B,A) \cong [B,A]_{\mathbf{C}}$. The reason for making
such fibrant-cofibrant replacements is that $\homs(B,A)$ is not a
priori homotopy invariant, in the sense that its weak homotopy
type may change if we replace $A$ or $B$ by a weakly equivalent
object. However, $\Mapbf(B,A)$ as defined above is homotopy
invariant. Of course, there are some choices involved in the
definition of $\Mapbf(B,A)$, but the weak homotopy type of
$\Mapbf(B,A)$ is well-defined.

In Section~\ref{S:2} we give an explicit model for
$\Mapbf(\mfH,\mfG)$ using weak 2-functors; see
Proposition~\ref{P:mappingspace2}.

\section{The pointed category $\mathbf{2Gpd}_*$ and  application
to crossed modules} {\label{S:Crossed1}}

In this section we discuss pointed versions of the results of the
previous section.
 The interest in the pointed version lies in the fact
that understanding the homotopy theory of {\em pointed} 2-groupoids
$\mathbf{2Gpd}_*$ is what we need in studying $\mathbf{2Gp}$ and
$\mathbf{CrossedMod}$.

  We define the  category $\mathbf{2Cat}_*$ of pointed 2-categories as follows.
  A pointed 2-category is a 2-category with a chosen base point (i.e.\ object)
  which,
   by abuse of notation, we denote by $*$.
  In $\mathbf{2Cat}_*$ a
  2-functor $P \: \mfD \to \mfC$ is required to preserve the base
  points, and a transformation $T \: P \Rightarrow Q$ between such
  functors is required to satisfy the condition $\theta_*=\id_*$
(similarly, for modifications we require $\mu_*$ to be the
identity). The pointed hom-2-category between two pointed
2-categories $\mfD$ and $\mfC$ is denoted by $\homs_*(\mfD,\mfC)$.

The same definitions apply to the category $\mathbf{2Gpd}_*$ of
2-groupoids. There is a pointed version of the Moerdijk-Svensson
closed model structure, and pointed versions of
Proposition~\ref{P:nerve} and Theorem~\ref{T:QuillenEq} are valid.

The following proposition follows from the pointed version of
Theorem~\ref{T:QuillenEq}.

\begin{prop}{\label{P:equiv}}
   The functor $N \: \Ho(\mathbf{2Gp}) \to \Ho(\mathbf{SSet}_*)$
   induces an equivalence between $\Ho(\mathbf{2Gp})$
   and the full subcategory of $\Ho(\mathbf{SSet}_*)$ consisting
   of connected pointed homotopy 2-types.
\end{prop}

We point out that every pointed 2-groupoid is fibrant in
$\mathbf{2Gpd}_*$. It is also easy to check that a pointed
2-groupoid is cofibrant in $\mathbf{2Gpd}_*$ if and only if it is
cofibrant in $\mathbf{2Gpd}$.

\subsection{Translation to the language of
crossed modules}{\label{SS:crossedmod}}

We can now apply the homotopy theory developed above for
$\mathbf{2Gpd}_*$ to study 2-groups. This is because $\mathbf{2Gp}$
naturally embeds as a full subcategory of $\mathbf{2Gpd}_*$, by
viewing a 2-group as a 2-groupoid with one object.  On the other
hand, 2-groups are the same as crossed modules (see
$\S$\ref{SS:Equiv}), so it is interesting to translate the
homotopical structure of $\mathbf{2Gpd}_*$ through $\mathbf{2Gp}$ to
the language of crossed modules.

We start with  the notions of transformation and pointed
transformation.
  Let $\mfG=[\varphi\: G_2 \to G_1]$ and $\mfH=[\psi \: H_2 \to
  H_1]$ be crossed modules, and let $P,Q \: \mfH \to \mfG$ be
  morphisms between them. Then a
    transformation  $T \: P \Rightarrow Q$ is given by
  a pair $(a,\theta)$ where $a \in G_1$ and $\theta \: H_1 \to G_2$
  is a crossed homomorphism  for the induced action, via $p_1^a$, of
  $H_1$ on $G_2$ (that is, $\theta(hh')=\theta(h)^{p_1(h')^a}\theta(h')$).
  These data should satisfy the following axioms:

    \begin{itemize}
        \item[$\mathbf{T1.}$] $p_1(h)^a\underline{\theta(h)}=q_1(h)$, for
             every $h \in H_1$;

        \item[$\mathbf{T2.}$] $p_2(\beta)^a\theta(\underline{\beta})=q_2(\beta)$,
             for every  $\beta \in H_2$.
    \end{itemize}
  A transformation $T$ is   pointed if and only if $a=1$.
  Between two pointed transformations
     there is no non-trivial {\em pointed} modification.

\begin{defn}
   We say a morphism $\mfH \to \mfG$ of crossed mod\-ules
   is a fibration (re\-spec\-tive\-ly, cofibration) if the
   induced morphism on the associated 2-groups is so
   (in the sense of Moerdijk-Svensson).
\end{defn}

\begin{prop}{\label{P:crossedfib}}
A map $(f_2,f_1) \: [H_2 \to H_1] \to  [G_2 \to G_1]$ of crossed
modules
   is a fibration if $f_2$ and $f_1$ are both surjective.
   It is  a trivial fibration if, furthermore, the
   map $H_2 \to H_1\times_{G_1} G_2$ is an isomorphism.
\end{prop}

Note that every crossed module is fibrant.

\begin{prop}{\label{P:crossedcofib}} A crossed module $[G_2 \to
G_1]$ is  cofibrant if and only if $G_1$
    is a free group.
\end{prop}

\begin{proof}
  This follows immediately from the Remark on page 194 of
  \cite{M-S}, but
  we give a direct proof.
  A 2-group $\mfG$ is cofibrant in the Moerdijk-Svensson structure
  if and only if every trivial fibration $\mfH \to \mfG$,
  where $\mfH$ is a 2-group{\em oid}, admits a section. But
  we can obviously restrict ourselves to 2-groups $\mfH$.
  So we can work entirely   within crossed modules, and
  use the crossed module version of trivial fibrations as in
  Proposition~\ref{P:crossedfib}.

   Assume $G_1$ is free. Let $(f_2,f_1)\: [H_2 \to H_1] \to  [G_2 \to G_1]$
   be a trivial fibration.
   Since $G_1$ is free and $f_1$ is surjective, there is a section
   $s_1 \: G_1 \to H_1$. Using  the fact that
   $H_2\cong H_1\times_{G_1} G_2$, we also get a natural section $s_2 \: G_2 \to H_2$
   for the projection $H_1\times_{G_1}G_2\to G_2$, namely,
   $s_2(\alpha)=(s_1\big(\underline{\alpha}),\alpha\big)$.
   It is easy to see that $(s_2,s_1) \: [G_2 \to G_1] \to [H_2 \to H_1]$
   is a map of crossed modules.

   To prove the converse, choose a free group $F_1$ and a surjection
    $f_1 \: F_1 \to G_1$.
    Form the pullback crossed module $[F_2 \to F_1]$
   by setting $F_2=F_1\times_{G_1} G_2$.  Then we have a trivial fibration
   $[F_2 \to F_1] \to  [G_2 \to G_1]$. By assumption, this has a section,
   so in particular we get a section $s_1 \: G_1 \to F_1$ which
   embeds $G_1$ as a subgroup of $F_1$. It follows from Nielsen's
   theorem that $G_1$ is free.
\end{proof}

\begin{rem}
   Let $G$ be a group and $H$ a subgroup. Let us say that $G$ is
   a {\em free extension} of $H$ if for every surjection
   $p \: K \to G$, every partial section $s \: H \to K$ to $p$ can be
   extended to $G$. For instance, if $H$ is the trivial group, this
   is equivalent to $G$ being free. Proposition~\ref{P:crossedfib}
   can be generalized by saying that
   $(f_2,f_1)\: [H_2 \to H_1] \to  [G_2 \to G_1]$
   is a cofibration if and only if $f_1 \: H_1 \to G_1$
   is injective and $G_1$ is a free extension of $f_1(H_1)$.
\end{rem}

Using the above proposition, we can  give a simple recipe for
cofibrant replacement of a crossed module. Let $\mfG=[G_2\to G_1]$
be an arbitrary crossed module.  Let $F_1 \to G_1$ be a surjective
map from a free group $F_1$,
      and set $F_2:=F_1\times_{G_1}G_2$. Consider the crossed module
      $\mfF=[F_2 \to F_1]$.  Then  $\mfF$ is cofibrant, and the natural
      map $\mfF \to \mfG$ is a trivial fibration (Proposition~\ref{P:crossedfib}).
      In other words, $\mfF \to \mfG$ is a cofibrant replacement for $\mfG$.

\section{Weak 2-functors between 2-categories}{\label{S:2}}

  In contrast to the case of ordinary groupoids, the nerve functor
  $N \: \mathbf{2Gpd} \to \mathbf{SSet}$ fails to be
  full. In order to make it full, one can cheat and define
  a {\em weak} 2-functor $\mfH \to \mfG$
  to be a simplicial map $N\mfH \to N\mfG$! Translating this
  back to the language of 2-categories, we arrive at
  the following definition.

  \begin{defn}{\label{D:weakmap1}}
    Let  $\mfH$ and $\mfG$ be 2-categories. A {\em weak}\/ 2-functor
    $F \:\mfH \to \mfG$ is defined the same way as a strict 2-functor,
    except that
    for every pair of composable arrows $a,b \in \mfH$ we require
    that instead of the equality $F(ab)=F(a)F(b)$, we are given
     a 2-morphism $\varepsilon_{a,b} \: F(a)F(b) \Rightarrow
     F(ab)$ which is
    natural in $a$ and $b$. For every three composable arrows $a$, $b$ and
    $c$, we require that the following diagram commutes:
     \[\xymatrix@R=12pt@C=32pt@M=10pt{
         F(a)F(b)F(c)\ar@{=>}[r]^{1_{F(a)}\cdot\varepsilon_{b,c}}
                                  \ar@{=>}[d]_{\varepsilon_{a,b}\cdot 1_{F(c)}}  &
                             F(a)F(bc)\ar@{=>}[d]^{\varepsilon_{a,bc}}  \\
          F(ab)F(c)   \ar@{=>}[r]_{\varepsilon_{ab,c}}        &    F(abc)   }\]
  \end{defn}

 \begin{rem}\label{R:weak}\indent\par
  \begin{itemize}
    \item[$\mathbf{i.}$] The naturalness assumption on
    $\varepsilon_{a,b}$ in particular implies the following.
    Suppose we are given  2-morphisms as in the diagram
    \[\xymatrix@C=18pt@R=2pt@M=2pt{ & \ar@{=>}[dd]^{\alpha}& & \ar@{=>}[dd]^{\beta}&\\
       \bullet  \ar@/^1.1pc/ [rr]^{a} \ar@/_1.1pc/ [rr]_{a'} &&
       \bullet  \ar@/^1.1pc/ [rr]^{b} \ar@/_1.1pc/ [rr]_{b'} && \bullet  \\
            & & & &    }\]
     Then, we require that after applying $F$ the diagram remains
     commutative. This means
         \[[\varepsilon_{a,b}][F(\alpha,\beta)]=[F(\alpha)F(\beta)][\varepsilon_{a',b'}],\]
      as  2-morphisms from $F(a)F(b)$ to $F(a'b')$.

    \item[$\mathbf{ii.}$] Our definition of weak is slightly stronger than the usual one in that
      we are assuming that  $F$ preserves the identity
      1-morphisms strictly. We do not, however, require that $F(a^{-1})=F(a)^{-1}$.
  \end{itemize}
\end{rem}

 Considering 2-groups as  2-categories with one object, we obtain
 the notion of a weak map of 2-groups. The translation
 of this definition in the language of crossed modules is given in
 Section~\ref{S:Crossed2}.

We denote the category of 2-categories (respectively, 2-groupoids,
2-groups) with weak morphisms by $\mathbf{2Cat}_{\text{w}}$
(respectively, $\mathbf{2Gpd}_{\text{w}}$,
$\mathbf{2Gp}_{\text{w}}$). Notice that each of these categories
contains the corresponding strict category as a subcategory.

\subsection{Transformations between weak
functors}{\label{SS:Transweak}}

Let $\mfD$ and $\mfC$ be 2-categories. Weak 2-functors from $\mfD$
to $\mfC$ form a 2-category. The 1-morphism in this 2-category are
transformations $T \: P \Rightarrow Q$. Recall
($\S$\ref{SS:Transformations2}) that a {\em transformation}
between strict functors $P$ and $Q$ is
 a rule which assigns to every morphism  $c \: A\to B$ in $\mfD$
 a 2-cell in $\mfC$ of the form
 \[\xymatrix@=12pt@M=10pt{
         P(A) \ar[r]^{t_A}\ar[d]_{P(c)}  &
              Q(A)  \ar[d]^{Q(c)}   \\
         P(B)    \ar[r]_{t_{B}}   \ar@{=>}[ur]^{\theta_c}     &  Q(B)   }\]
In case of weak 2-functors $P$ and $Q$, we require the same
conditions as in $\S$\ref{SS:Transformations2}, except that the
second compatibility condition ``$\theta_a\theta_b=\theta_{ab}$'',
should be modified as follows.
 For every two
composable morphisms $a$ and $b$, the prism (of the form
$\Delta^1\times\Delta^2$) with three square faces  $\theta_a$,
$\theta_b$ and $\theta_{ab}$ and two triangular faces
$\epsilon_{a,b}^P$ and $\epsilon_{a,b}^Q$ is commutative.

We can also talk about modification between  two weak
2-transformations between weak 2-functors. The definition  is
identical to the one in the case of strict 2-functors
($\S$\ref{SS:Transformations}).

\begin{defn}{\label{D:weakfunctors}}
  Given 2-categories $\mathfrak{C}$ and $\mathfrak{D}$,
  we define $\homw(\mathfrak{D},\mathfrak{C})$ to be the
  2-category whose objects are weak 2-functors from
  $\mathfrak{D}$ to $\mathfrak{C}$, whose 1-morphisms
  are weak 2-transformations between 2-functors, and whose
  2-morphisms are modifications.
\end{defn}

Observe that $\homw(\mathfrak{D},\mathfrak{C})$ contains
$\homs(\mathfrak{D},\mathfrak{C})$ as a full sub-2-category. In the
case where $\mathfrak{C}$ is a 2-groupoid,
$\homw(\mathfrak{D},\mathfrak{C})$ is also a 2-groupoid.

\subsection{Nerves of weak functors}

As pointed out at the beginning of this section, our definition of a
weak functor between two 2-categories was formulated by translating
the simplicial identities that have to be satisfied by a simplicial
map on the corresponding nerves. Having taken this for granted, the
following is immediate.

\begin{prop}{\label{P:weak}}
   The nerve functors extend to  fully faithful
  functors \[N \: \mathbf{2Cat}_{\text{w}} \to \mathbf{SSet}, \
  N \: \mathbf{2Gpd}_{\text{w}} \to \mathbf{SSet}, \  \text{and} \
  N \: \mathbf{2Gp}_{\text{w}} \to \mathbf{SSet}_*.\]
\end{prop}

\begin{rem}{\label{R:nerve}}
 Using the notion of  weak functor between 2-categories, we
 can give an alternative view to the nerve of a 2-category $\mfG$.
 Let $\mathbf{I}_n$ be the or\-dered set $\{0,1,\cdots,n\}$, with the usual
 or\-der\-ing, viewed as a cat\-e\-gory (hence, also a 2-cat\-e\-gory). Then
 $N\mfG_n$ is in natural bijection with the set of weak functors
 $\mathbf{I}_n \to \mfG$. This also explains why $N$ is functorial
 and  preserves fiber products.
\end{rem}

\begin{lem}{\label{L:weakweak}}
  A weak map $P \: \mfH \to \mfG$ of  2-groupoids induces
  group homomorphisms  $\pi_iP\: \pi_i\mfH \to \pi_i\mfG$, $i\leq 2$, for every choice
  of base points. In particular,
  we can talk about a weak map $P$ being an equivalence of
  2-groupoids. (Compare Definition~\ref{D:homotopygroup}.)
  Furthermore, $P \: \mfH \to \mfG$ is an equivalence of 2-groupoids
  if and only if $NP \: N\mfH \to N\mfG$ is a weak equivalence of
  simplicial sets.
\end{lem}
\begin{proof}
 Straightforward.
\end{proof}

The  next result is that $\homw(\mfH,\mfG)$ can be used to compute
derived hom-2-groupoids without the need to take a cofibrant
replacement for $\mfH$. To prove this fact, first we prove a lemma.

\begin{lem}{\label{L:homotopy}}
 Let $P,Q \: \mfH\to\mfG$ be weak 2-functors between 2-groupoids. Then $P$ and
 $Q$,
 viewed as objects in the 2-groupoid $\homw(\mfH,\mfG)$, are isomorphic if and only if
 $NP$ and $NQ$ are simplicially homotopic. In other words, we have
 a natural bijection
   \[\pi_0\homw(\mfH,\mfG) \cong
   [N\mfH,N\mfG]_{\mathbf{SSet}}.\]
\end{lem}

\begin{proof}
    Observe that a simplicial homotopy from $NP$ to $NQ$,
    that is, a simplicial map $\Delta^1\times N\mfH \to \mfG$
    connecting $NP$ to $NQ$, is given by data almost identical
    to that required to give a weak transformation
    from $P$ to $Q$ (see $\S$\ref{SS:Transformations2}), where
    instead of the squares of $\S$\ref{SS:Transformations2}
    we have to use squares of the form
    \[\xymatrix@=6pt@M=10pt{
         P(A) \ar[rr]^{t_A}\ar[dd]_{P(c)}  \ar[rrdd] |-{t_c} &&
        Q(A)  \ar[dd]^{Q(c)}   \ar@{=>}[dl]^{\uc{\theta}{c}}\\ & & \\
         P(B)    \ar[rr]_{t_{B}}   \ar@{=>}[ur]^{\lc{\theta}{c}}
        & &  Q(B)   }\]
  Namely, the 2-morphism $\theta_c$ is now replaced by a triple
  $(t_c,\lc{\theta}{c},\uc{\theta}{c})$. The compatibility
  conditions required for $(t_c,\lc{\theta}{c},\uc{\theta}{c})$
  are simply translated, in the obvious way,  from the ones for
  $\theta_c$.

  It is now easy to prove the lemma. Given a simplicial homotopy
  from $NP$ to $NQ$, presented by  triples
  $(t_c,\lc{\theta}{c},\uc{\theta}{c})$,
  one defines
  $\theta_c:=[\lc{\theta}{c}][\uc{\theta}{c}]^{-1}$. This gives a
  transformation from $P$ to $Q$. Conversely, given a
  transformation from $P$ to $Q$, presented by $\theta_c$, one
  defines the triples $(t_AQ(c),\theta_c,\id)$. This
  gives a simplicial homotopy from $NP$ to $NQ$.
\end{proof}

\begin{prop}{\label{P:mappingspace2}}
 Let $\mfG$ and $\mfH$ be 2-groupoids. Then there is a
 natural (up to weak transformation)   equivalence of 2-groupoids
      \[\Mapbf(\mfH,\mfG) \simeq \homw(\mfH,\mfG).\]
 In other words, $\homw(\mfH,\mfG)$ is a natural model
 for the derived hom-2-groupoid of functors from $\mfH$ to $\mfG$.
 In particular, when $\mfH$ is cofibrant, we have a natural
 equivalence of 2-groupoids
      \[\homs(\mfH,\mfG) \simeq \homw(\mfH,\mfG).\]
\end{prop}

\begin{proof}
   First we show that if $p \: \mfH' \to \mfH$ is an arbitrary
  equivalence, then the induced map
   $p^* \: \homw(\mfH,\mfG) \to \homw(\mfH',\mfG)$ is also an
   equivalence. Since the induced map $Np \: N\mfH' \to N\mfH$ is
   a homotopy equivalence of Kan complexes,
   by Whitehead's theorem, it admits an inverse
   $Q \: N\mfH \to N\mfH'$, in the sense that $QNp$ and
   $Np\,Q$ are homotopic to the corresponding identity maps.
   Since $N$ is fully faithful
   (Proposition~\ref{P:weak}), there
   exists a weak functor $q \: \mfH \to \mfH'$ such that
   $Nq=Q$. It follows from Lemma~\ref{L:homotopy} that $p$ and $q$
   are inverse equivalences. That is, each of
   $f\circ q \: \mfH \to \mfH$ and $q\circ p \: \mfH' \to \mfH'$
   is connected to the corresponding identity functor by a
   weak transformation. This implies that
   $p^* \: \homw(\mfH,\mfG) \to \homw(\mfH',\mfG)$
   and $q^* \: \homw(\mfH',\mfG) \to \homw(\mfH,\mfG)$ are also
   inverse equivalences.

   To prove the proposition, we may assume now that
   $\mfH$ is cofibrant. We have a
   natural full  inclusion $\homs(\mfH,\mfG) \subset
   \homw(\mfH,\mfG)$. So it is enough to show that
   the induced map $\pi_0\homs(\mfH,\mfG) \to
   \pi_0\homw(\mfH,\mfG)$ is a bijection. Since $\tgpd$
   is a monoidal closed model category and $\mfH$ is cofibrant,
   we have that
    $\pi_0\homs(\mfH,\mfG)\cong [\mfH,\mfG]_{\mathbf{2Gpd}}$.
   On the other hand, by Lemma~\ref{L:homotopy}, we have that
    \[\pi_0\homw(\mfH,\mfG) \cong
   [N\mfH,N\mfG]_{\mathbf{SSet}}.\]
   The claim now follows
   from the fact that $N \: \Ho(\tgpd) \to \Ho(\sset)$ is an equivalence
   of categories (which, in particular, implies that $N$ induces a bijection between
   $[\mfH,\mfG]_{\mathbf{2Gpd}}$ and
   $[N\mfH,N\mfG]_{\mathbf{SSet}}$).
\end{proof}

\begin{prop}{\label{P:mappingspace3}}
    Given 2-groupoids $\mfG$ and $\mfH$, there is a natural
  homotopy equivalence of simplicial sets
    \[N\Mapbf(\mfH,\mfG) \simeq \Hombf(N\mfH,N\mfG).\]
 \end{prop}

\begin{proof}
 Making a cofibrant replacement for $\mfH$ does not change either
 side, up to a canonical homotopy, so we
 may assume $\mfH$ is cofibrant.
 We have to show that
     \[N\homs(\mfH,\mfG) \simeq \Hombf(N\mfH,N\mfG).\]
 We prove that for every simplicial set $A$, there is a natural
 bijection
   \[[A,N\homs(\mfH,\mfG)]_{\mathbf{SSet}} \cong
     [A,\Hombf(N\mfH,N\mfG)]_{\mathbf{SSet}}.\]
 That is, the contravariant functors from $\mathbf{SSet}$ to $\mathbf{Set}$
 represented by $N\homs(\mfH,\mfG)$ and  $\Hombf(N\mfH,N\mfG)$
 are naturally isomorphic.

 Since both $N\homs(\mfH,\mfG)$ and $\Hombf(N\mfH,N\mfG)$ have
 trivial homotopy group in degrees   higher than 2, we may assume
 that $A$ is a simplicial set with the same property. Since the unit of
 adjunction $A \to NW(A)$ is a weak equivalence
 (Theorem~\ref{T:QuillenEq}.$\mathbf{iii}$), we may assume
 $A=N\mfK$, for some cofibrant 2-groupoid $\mfK$. We start
 by simplifying  the right hand side:
 \begin{eqnarray}
   [A, N\homs(\mfH,\mfG)]_{\mathbf{SSet}}
      & = & [N\mfK, N\homs(\mfH,\mfG)]_{\mathbf{SSet}}  \nonumber  \\
        & \cong & [\mfK, \homs(\mfH,\mfG)]_{\mathbf{2Gpd}}  \nonumber\\
       & \cong & \pi_0\homs\big(\mfK,\homs(\mfH,\mfG)\big)    \nonumber\\
       & \cong & \pi_0\homs(\mfK\otimes \mfH,\mfG)  \nonumber\\
       & \cong & [\mfK\otimes \mfH,\mfG]_{\mathbf{2Gpd}}  \nonumber\\
       & \cong & [\mfK\times \mfH,\mfG]_{\mathbf{2Gpd}}.  \nonumber
 \end{eqnarray}
 The next to last equality follows from the fact that in every
 monoidal model category, the tensor product of two cofibrant objects
 is cofibrant.
 The last equality follows from the fact that, for arbitrary
 2-groupoids
 $\mfK$ and $\mfH$, the natural morphism $\mfK\otimes \mfH \to \mfK\times
 \mfH$ is an equivalence of 2-groupoids (in fact it is a trivial
 fibration); see the last paragraph of (\cite{Lack}, Section 2).

Now we do the left hand side:
 \begin{eqnarray}
      [A,\Hombf(N\mfH,N\mfG)]_{\mathbf{SSet}}
         & = & [N\mfK,\Hombf(N\mfH,N\mfG)]_{\mathbf{SSet}} \nonumber  \\
         & \cong &  \pi_0\big(\Hombf(N\mfK,\Hombf(N\mfH,N\mfG)\big) \nonumber  \\
         & \cong &  \pi_0\big(\Hombf(N\mfK\times N\mfH,N\mfG)\big)   \nonumber  \\
         & \cong &  [N\mfK\times N\mfH,N\mfG]_{\mathbf{SSet}} \nonumber  \\
         & \cong &  [N(\mfK\times \mfH),N\mfG]_{\mathbf{SSet}} \nonumber  \\
         & \cong &  [\mfK\times \mfH,\mfG]_{\mathbf{2Gpd}}. \nonumber
 \end{eqnarray}

We see that the left hand side coincides with the right hand side.
This implies that $N\homs(\mfH,\mfG)$ and $\Hombf(N\mfH,N\mfG)$
are naturally isomorphic in $\Ho(\mathbf{SSet})$. Since both
 $N\homs(\mfH,\mfG)$ and $\Hombf(N\mfH,N\mfG)$ are Kan complexes,
 it follows from Whitehead's theorem that they are naturally homotopy
 equivalent.
\end{proof}

\begin{cor}{\label{C:mappingspace}}
   Let $\mfG$ and $\mfH$ be 2-groupoids. Then, there is a
   natural (up to homotopy) homotopy equivalence of simplicial sets
       \[\Hombf(N\mfH,N\mfG)\simeq N\Mapbf(\mfH,\mfG)
                             \simeq N\homw(\mfH,\mfG).\]
\end{cor}

\begin{proof}
  Combine Proposition~\ref{P:mappingspace2} and
  Proposition~\ref{P:mappingspace3}.
\end{proof}

\begin{cor}{\label{C:lift}}
  Let $\mfG$ and $\mfH$ be 2-groupoids, and let
  $f \in [\mfH,\mfG]_{\mathbf{2Gpd}}$.
  Then there is a weak map
  $\tilde{f} \: \mfH \to \mfG$,
  unique up to transformation, which induces $f$ in the homotopy
  category.
\end{cor}

\begin{proof}
  This follows from the fact that
    \[[\mfH,\mfG]_{\mathbf{2Gpd}} \cong
    \pi_0\Mapbf(\mfH,\mfG)\simeq\pi_0\homw(\mfH,\mfG).\]
\end{proof}

In (\cite{M-S}, Proposition 2.2.i) it is  claimed that the functor
$W \: \mathbf{SSet} \to \mathbf{2Gpd}$, the left adjoint to the
nerve functor, preserves products. Unfortunately this is easily
seen to be false: $W(I\times I)$ is not isomorphic to $W(I)\times
W(I)$ (but they are naturally homotopy equivalent). If this were
true, an easy   argument would imply that there is natural
isomorphism of simplicial sets
     \[N\homs(WX,\mfG) \cong \Hombf(X,N\mfG),\]
for every  simplicial set $X$ and every  2-groupoid $\mfG$; in
other words, the adjunction between $W$ and $N$ would be an
enriched adjunction (which is not true either). Since $WX$ is
cofibrant (\cite{M-S}, Proposition 2.2.ii), we would indeed have
    \[N\Mapbf(WX,\mfG) \cong N\homs(WX,\mfG) \cong \Hombf(X,N\mfG),\]
which would lead to a straightforward proof of Proposition
\ref{P:mappingspace3}.\footnote{There is a similar type of error
in \oldcite{M-S} in the Remark at the bottom of page 194 where the
authors claim that the interval groupoid $\mathfrak{I}$ yields a
cylinder object $\mathfrak{I} \times A$ for every 2-groupoid $A$.
This is false, because the map $A\coprod A \to \mathfrak{I} \times
A$ is not in general a cofibration. The right thing to do would be
to take $\mathfrak{I} \otimes A$. Or one could work with  path
objects $\mfG^{\mathfrak{I}}:=\homs(\mathfrak{I},\mfG)$.}

This line of argument is unfortunately flawed.  Nevertheless, one
can prove the following result (which will not be needed anywhere
in this paper).

\begin{prop}{\label{P:mappingspace4}}
    Given a 2-groupoid $\mfG$ and a simplicial set $X$, there is a
    natural homotopy equivalence of simplicial sets
       \[N\Mapbf(WX,\mfG) \simeq \Hombf(X,N\mfG).\]
    That is, the adjunction between $N$ and $W$ can be enriched,
    up to homotopy.
\end{prop}

\begin{proof}
  Note that since $\mfG$ has trivial homotopy groups in degrees
  higher than 2, and since the unit of adjunction $X \to NW(X)$
  induces isomorphisms on $\pi_i$, $i\leq 2$, the two simplicial sets
  $\Hombf(X,N\mfG)$ and $\Hombf(NW(X),N\mfG)$ are naturally
  homotopy  equivalent. Now apply Proposition~\ref{P:mappingspace3} with
  $\mfH=W(X)$.
\end{proof}

\subsection{Monoidal properties of $N$ and $W$}{\label{SS:LaxMon}}

In this subsection we say a few words about the monoidal properties
of the functors $N$ and $W$. We briefly touch upon this issue
without getting into much detail. The material in this subsection
will not be used elsewhere in the paper.

What follows is for the most part a consequence of Tonks'
Eilenberg-Zilber theory \cite{Tonk}. Tonks' results are, however,
stated in the language of crossed complexes of groupoids (while we
work with 2-groupoids). Consequently, his monoidal structure  is
given by the Brown-Higgins tensor product of crossed complexes
(while we use the Gray tensor product of 2-groupoids). To make the
passage from the category $\mathbf{Crs}$ of crossed modules in
groupoids to the category $\mathbf{2Gpd}$ of 2-groupoids, we observe
that the latter can be identified with a full subcategory of the
former. Furthermore, there is an idempotent ``truncation functor''
$t \: \mathbf{Crs} \to \bf{2Gpd}$ which sends the Brown-Higgins
tensor to the Gray tensor, namely  $t(A\otimes_{BH}B)=t(A)\otimes
t(B)$. Using the functor $t$ we can translate Tonks' results to the
language of 2-groupoids.

\medskip

Let us now outline certain monoidal features of the functors $N$ and
$W$. The functors $N$ and $W$ are both lax monoidal and lax
comonoidal, and these lax monidal/comonoidal structures on these
functors are compatible in many ways:

\medskip

\noindent {\bf $W$ is lax monoidal.} There are  morphisms of
2-groupoids $b_{X,Y} \: W(X)\otimes W(Y) \to W(X\times Y)$, natural
in simplicial sets $X$ and $Y$. This follows from \cite{Tonk}, Proposition 2.6.

\medskip

\noindent {\bf $W$ is lax comonoidal.} There are  morphisms of
2-groupoids $a_{X,Y} \:  W(X\times Y) \to W(X)\otimes W(Y)$, natural
in simplicial sets $X$ and $Y$. This follows from \cite{Tonk}, Proposition 2.1.

\medskip

\noindent {\bf $N$ is lax monoidal.} There are  morphisms of
simplicial sets $a'_{G,H} \: N(G)\times N(H) \to N(G\otimes H)$,
natural in $G$ and $H$. This  follows from \cite{Tonk}, Corollary
2.4.

\medskip

\noindent {\bf $N$ is lax comonoidal.} There are  morphisms of
simplicial sets $b'_{G,H} \:  N(G\otimes H) \to  N(G)\times N(H)$,
natural in $G$ and $H$. To get these morphisms, simply apply $N$ to
the natural morphism $G \otimes H \to G\times H$, and use the fact
that $N(G\times H)=NG \times NH$.

\medskip

\noindent {\bf Compatiblity.}  The monoidal and comonoidal
structures on $W$ are related by (a straightforward modification of)
the Eilenberg-Zilber theorem of Tonks (\cite{Tonk}, Theorem 3.1):
the compositions $a_{X,Y} \circ b_{X,Y}$ are equal to the identity,
and the compositions $b_{X,Y}\circ a_{X,Y}$ are homotopic to
identity---homotopy meaning left-homotopy, i.e.,  via mapping
cylinders $G \sst{} I\otimes G$---through homotopies $\phi_{X,Y}$
that are natural in $X$ and $Y$ and fix the image of $b$
 \[\xymatrix@=16pt@M=8pt@C=28pt{ \ar@(ul,dl)[]|{\phi} }
      \xymatrix@=16pt@M=8pt@C=28pt{  W(X\times Y)  \ar@<1ex> [r]^a
       & \ar@<1ex> [l]^b WX\otimes WY. }\]
The same thing can be said about $N$. In other words, $W(X)\otimes
W(Y)$ is naturally a strong deformation retract of $W(X\times Y)$,
and $N(G)\times N(H)$ is naturally a strong deformation retract of
$N(G\otimes H)$. In particular, all morphism $a$, $b$, $a'$ and $b'$
are weak equivalences.

It also seems plausible that the pair $(W,N)$ is indeed a {\em
monoidal adjunction} (and also a {\em comonoidal adjunction}). One
could also ask whether the monoidal (resp., comonoidal) structure of
$W$ is compatible with the comonoidal (resp., monoidal) structure of
$N$. One may also wonder if the above mentioned deformation
retractions for $W$ are in some way compatible with the ones for
$N$.


\section{Pointed weak functors and application to crossed modules}{\label{S:Crossed2}}

We can prove pointed versions of the results of the previous
section. For instance, we have the following propositions.

\begin{prop}{\label{P:pointedmappingspace2}}
 Let $\mfG$ and $\mfH$ be pointed 2-groupoids. Then, there is a
 natural (up to homotopy)   equivalence of 2-groupoids
      \[\Mapbf_*(\mfH,\mfG) \simeq \homw_*(\mfH,\mfG).\]
 In other words, $\homw_*(\mfH,\mfG)$ is a natural model
 for the derived groupoid of pointed  functors from $\mfH$ to $\mfG$.
 In particular, when $\mfH$ is cofibrant, we have a natural
 equivalence of 2-groupoids
      \[\homs_*(\mfH,\mfG) \simeq \homw_*(\mfH,\mfG).\]
\end{prop}

\begin{prop}
   Let $\mfG$ and $\mfH$ be pointed 2-groupoids. Then, there is a
   natural (up to homotopy) homotopy equivalence of simplicial sets
       \[\Hombf_*(N\mfH,N\mfG)\simeq N\Mapbf_*(\mfH,\mfG) \simeq N\homw_*(\mfH,\mfG).\]
\end{prop}

\begin{proof}
    We prove that whenever $\mfH$ is a cofibrant pointed 2-groupoid and
   $\mfG$ an arbitrary pointed 2-groupoid, then we have a natural
   homotopy equivalence
     \[N\homs_*(\mfH,\mfG) \simeq \Hombf_*(N\mfH,N\mfG).\]
   We will always denote the base points by $*$.
   Consider the map of 2-groupoids
     \[\homs(\mfH,\mfG) \llra{\operatorname{ev}_*} \mfG.\]
   Applying Axiom $\mathbf{M7}$ of monoidal closed model categories
   (\cite{Lack}, 7)   to the cofibration $* \to \mfH$
   and fibration $\mfG \to *$, we see that
   $\operatorname{ev}_*$ is a fibration of 2-groupoids.
   Similarly, using the fact that the nerve of a 2-groupoid is fibrant
   (Proposition~\ref{P:nerve}.$\mathbf{ii}$), together
   with $\mathbf{SM7}$ for $\mathbf{SSet}$,
   we see that the following map of simplicial sets is also
   a fibration:
     \[\Hombf(N\mfH,N\mfG) \llra{\operatorname{ev}_*} N\mfG.\]
   Now, we observe that $\homs_*(\mfH,\mfG)$ is the fiber over
   $*$ of $\operatorname{ev}_*\: \homs(\mfH,\mfG) \to \mfG$,
   and $\Hombf_*(N\mfH,N\mfG)$ is
   the fiber over $*$ of $\operatorname{ev}_* \: \Hombf(N\mfH,N\mfG) \to
   N\mfG$. The claim now follows from
   Proposition~\ref{P:mappingspace3}, plus the fact that taking nerves preserves
   weak equivalences and fiber products
   (Proposition~\ref{P:nerve}.$\mathbf{i}$,$\mathbf{ii}$).
\end{proof}

\begin{prop}{\label{P:pointedlift}}
  Let $\mfG$ and $\mfH$ be pointed 2-groupoids. Then, for every
  $f \in\break [\mfH,\mfG]_{\mathbf{2Gpd}_*}$
  there is a weak pointed map $\tilde{f} \: \mfH \to \mfG$,
  unique up to pointed transformation, which induces $f$ in the homotopy
  category.
\end{prop}

Using the fact that $\mathbf{2Gp}$ is a full subcategory of
$\mathbf{2Gpd}_*$, we immediately obtain the 2-group versions of the
above propositions. In particular, using the fact that
$\mathbf{CrossedMod}$ and $\mathbf{2Gp}$ are equivalent
($\S$\ref{SS:Equiv}), we get the corresponding statements for
crossed modules as well. Note that in this case both sides of the
equivalence in Proposition~\ref{P:pointedmappingspace2} are actually
1-groupoids.

For the readers' information we include the crossed module
translations of the notions of  weak functor between 2-group,
transformation between weak functors, and modification between
transformations.

  \begin{defn}{\label{D:crossedweakmap}}

  Let $\mfG$ and $\mfH$ be crossed modules.
    A {\em weak map} $P \: \mfH \to \mfG$ consists of the following data:
   \begin{itemize}
   \item a pointed set map $p_1 \: H_1 \to G_1$,

   \item a pointed set map $p_2 \: H_2 \to G_2$,

   \item a set map $\varepsilon \: H_1\times H_1 \to G_2$, denoted
         by $(x,y) \mapsto \varepsilon_{x,y}$.
\end{itemize}
   These data should satisfy the following conditions:
   \begin{itemize}
   \item[$\mathbf{W1.}$] \ $\forall \alpha \in H_2$,
              $p_1(\underline{\alpha})=\underline{p_2(\alpha)}$;

   \item[$\mathbf{W2.}$] \ $\forall \alpha, \beta\in H_2$,
      $p_2(\alpha\beta)=p_2(\alpha)p_2(\beta)
          \varepsilon_{\underline{\alpha},\underline{\beta}}$;

   \item[$\mathbf{W3.}$] \ $\forall x,y \in H_1$,
                 $p_1(xy)=p_1(x)p_1(y)\varepsilon_{x,y}$;

   \item[$\mathbf{W4.}$] Cocycle condition:
     $$\forall x,y,z \in H_1,\ \
            \varepsilon_{x,y}^{p_1(z)}\varepsilon_{xy,z}=
                       \varepsilon_{y,z}\varepsilon_{x,yz};$$

   \item[$\mathbf{W5.}$]  Equivariance:
       $$\forall x\in H_1, \ \forall \beta\in H_2,\ \
        \varepsilon_{x^{-1},x}p_2(\beta^x)=
                      p_2(\beta)^{p_1(x)}\varepsilon_{\underline{\beta},x}
                       \varepsilon_{x^{-1},\underline{\beta}x}.$$

  \end{itemize}
\end{defn}

  Ettore Aldrovandi has found a two-variable version of
  ($\mathbf{W5}$) that, in the presence of the other axioms,
  is equivalent to ($\mathbf{W5}$).

  \begin{itemize}

    \item[$\mathbf{W5'.}$] Two-variable version:
       $$\forall x,y \in H_1, \ \forall \beta \in H_2,\ \
           \varepsilon_{y,x}p_2(\beta^x)\varepsilon_{yx,x^{-1}\underline{\beta}x}=
                      p_2(\beta)^{p_1(x)}\varepsilon_{\underline{\beta},x}
                                         \varepsilon_{y,\underline{\beta}x}.$$
  \end{itemize}

\begin{defn}{\label{D:crossedtransfo}}
Suppose we are given two weak functors $P=(p_1,p_2,\epsilon)$ and
$Q=(q_1,q_2,\delta)$ from $\mfH$ to $\mfG$, as in Definition
\ref{D:crossedweakmap}.
 A {\em transformation} from $P$ to $Q$ is
 given by a pair $(a,\theta)$, where $a \in G_1$, and
$\theta \: H_1 \to G_2$ is a map of sets. This pair should satisfy
the following conditions:

 \begin{itemize}
     \item[$\mathbf{T0.}$] \
     $\varepsilon_{x,y}^a\theta(xy)=\theta(x)^{p_1(y)^a}\theta(y)\delta_{x,y}$,
        for every $x,y \in H_1$;

     \item[$\mathbf{T1.}$] \ $p_1(x)^a\underline{\theta(x)}=q_1(x)$, for
             every $x \in H_1$;

     \item[$\mathbf{T2.}$] \ $p_2(\alpha)^a\theta(\underline{\alpha})=q_2(\alpha)$,
             for every  $\alpha \in H_2$.
 \end{itemize}
 Such a transformation is called {\em pointed} if $a=1$.
\end{defn}

\begin{defn}{\label{D:crossedmod}}
  A {\em modification} between two transformations $(a,\theta)$
  and $(b,\sigma)$ is given by an element $\mu \in G_2$ that has
  the
  following properties:
    \begin{itemize}
     \item[$\mathbf{M1.}$] \ $a\underline{\mu}=b$;

     \item[$\mathbf{M2.}$] \ $\mu\sigma(x)=\theta(x)\mu^{q_1(x)}$, for every $x \in H_1$.
    \end{itemize}
  By definition, the only pointed modification is $\mu=1$, namely
      the trivial modification.

\end{defn}

Crossed modules are in practice more amenable to computations than
2-groups, while 2-groups are more conceptual. An important problem
regarding 2-groups is to compute the derived hom 2-groupoids
$\Mapbf(\mfH,\mfG)$ in an explicit way. One way to do this is to use
the above description of weak morphisms. It is also possible to give
a description of elements of $\Mapbf(\mfH,\mfG)$ using cohomological
invariants. But this may in practice be not  so useful, and we also
loose information about transformations between weak morphisms. In
\cite{Maps} we give an explicit description of weak morphisms, and
transformations between them, in terms of certain diagrams of groups
which we call {\em butterflies}.

\section{Weak  2-groupoids}{\label{S:3}}

In Section~\ref{S:1} we saw (Proposition~\ref{P:nerve}) that the
nerve functor $N \: \mathbf{2Gpd} \to \mathbf{SSet}$ identifies
$\mathbf{2Gpd}$ with a subcategory of $\mathbf{SSet}$ that is not
full. The lack of fullness was remedied in Section~\ref{S:2} by
considering weak functors between 2-groupoids (Proposition
\ref{P:weak}). So the category $\mathbf{2Gpd}_w$ is identified with
a certain full subcategory of $\mathbf{SSet}$. The natural question
now is to characterize this full subcategory. We know\pagebreak\ so far that a
simplicial set $X$ in the image of $N$ has the following properties:
 \begin{itemize}
   \item[$\mathbf{1.}$] $X$ is a Kan complex.

   \item[$\mathbf{2.}$] $X$ is 3-coskeletal (i.e., $X$ is isomorphic
              to the coskeleton of its 3-truncation $\{X_0,X_1,X_2,X_3\}$.
  \end{itemize}

  Is this enough to characterize the image of $N$? The answer is no.
  To see what is missing, let us consider the nerve
  functor from $\mathbf{Gpd}$ to $\mathbf{SSet}$.
  In this case, every simplicial set in
  the image of $N$ is
  Kan and 2-coskeletal. But a simplicial set in the image of $N$ has
  an additional property:
  it satisfies the condition for a {\em minimal} simplicial set for $n \geq 1$.
  That is, such $X$ is  minimal, except possibly at degree 0.

  \begin{defn}{\label{D:minimal}}
    We say a simplicial set $X$ is $k$-{\em minimal} if it satisfies the
    minimality condition for $n \geq k$. That is, every two $n$-simplices
    that have  equal  boundaries and  are homotopic relative to
    their common boundary are actually equal
    (see\cite{May},  \S 9 , or \cite{Jardine}, \S I.10).
  \end{defn}

  For example, the nerve  of a groupoid is 1-minimal, as we explained above.

  It is easy to see that any $X$ that is isomorphic to the nerve of a 2-groupoid
  satisfies the following property:

  \begin{itemize}
    \item[$\mathbf{3.}$] $X$ is 2-minimal.
  \end{itemize}

   \begin{rem}\label{R:minimal}\indent\par
   \begin{itemize}
   \item[$\mathbf{a.}$] Any $k$-coskeletal $X$ is
      automatically ($k+1$)-minimal.
      So, any
      3-skeletal $X$ is 4-minimal.
      Therefore, the condition ($\mathbf{3}$)
      puts only two new restrictions on $X$. In other words, in the presence of
      ($\mathbf{1}$) and ($\mathbf{2}$), we only need
      to verify that the condition of Definition~\ref{D:minimal}
      is satisfied with $n=2$ and $n=3$.

   \item [$\mathbf{b.}$]If $X$ is Kan and $2$-minimal, then it can be shown that
      for any $n\geq 3$,
      whenever two $n$-simplices in $X$ have the same boundary,
      then they are indeed
      homotopic relative boundary (hence equal).
      In other words, the natural map $X \to \Cosk_2(X)$ is injective
      (in general, $k$-minimal implies $X \to \Cosk_kX$ injective).
      So, one of the conditions in ($\mathbf{a}$) can be replaced
      by the injectivity of $X \to \Cosk_2X$.
  \end{itemize}
  \end{rem}

  Let $\mathbf{SSet2}$ denote the full-subcategory of $\mathbf{SSet}$
  consisting of simplicial sets satisfying  conditions
  ($\mathbf{1}$), ($\mathbf{2}$) and ($\mathbf{3}$)
  above. It is more natural now to consider
  the (fully faithful) functor
  $N \: \mathbf{2Gpd}_w \to \mathbf{SSet2}$ and ask whether
  it is an equivalence. The answer is no! In the same way we had to weaken the
  notion of a functor in order to make $N$ full, we have to weaken the notion
  of a 2-groupoid to make $N$ essentially surjective.

\begin{defn}{\label{D:weak2gpd}}
   A {\em weak 2-category} $\mfC$ is defined by exactly the same
   data and axioms as a strict 2-category, except that the axiom
   on the associativity of composition of three composable
   arrows $a$, $b$ and $c$ is replaced by the choice of an associator
   $\phi_{a,b,c} \: (ab)c \Rightarrow a(bc)$, for every such triples of
   arrows.\footnote{The properties
   of the identity arrows, however,
   remain strict. In this sense, our notion of weak is stronger
   than the usual one.}\pagebreak\ The associators are
   required to satisfy the following conditions:
   \begin{itemize}
      \item[$\mathbf{A1.}$] The famous pentagon identity
        \[[\phi_{a,b,c}d][\phi_{a,bc,d}][a\phi_{b,c,d}]
                          =[\phi_{ab,c,d}][\phi_{a,b,cd}].\]
       See, for example (\cite{B}, page 6).

      \item[$\mathbf{A2.}$] Whenever one of $a$, $b$ or $c$ is the
      identity, then $\phi_{a,b,c}$ is the identity.
   \end{itemize}
   A {\em weak 2-groupoid} is a
   weak 2-category in which for every arrow $a$ there exists
   an arrow $b$, with reverse source and target as $a$,
   such that $ab$ and $ba$ are   2-isomorphic to the corresponding identity arrows.
   A {\em weak 2-group} is a weak 2-groupoid with one object.
\end{defn}

\begin{rem}{\label{R:inverse}}
   The requirement on   inverses seems a little too weak, as we have not
   asked for any compatibility with  the associators.
   However, as is shown in \cite{Baez} for the case of weak 2-groups,
   it can be shown that one can {\em select} a coherent family of inverses in
   a way that all the required compatibility conditions with the associators
   are satisfied. We will not, however, need this fact here.
\end{rem}

\subsection{Weak functors between weak 2-categories, and
transformations}{\label{SS:weakfunct}}

We can  talk about {\em weak functors} between weak 2-categories.
The definition is the same as Definition~\ref{D:weakmap1}, but the
square should be replaced by the following hexagon
   \[\xymatrix@R=16pt@C=20pt@M=7pt{
          & F(a)\big(F(b)F(c)\big) \ar@{=>}@<1ex>[rd]^(0.6){1_{F(a)}\cdot\varepsilon_{b,c}} & \\
             \big(F(a)F(b)\big)F(c)\ar@{=>}@<1ex>[ru]^{\phi}
                         \ar@{=>}[d]_{\varepsilon_{a,b}\cdot 1_{F(c)}}
          &&       F(a)F(bc)\ar@{=>}[d]^{\varepsilon_{a,bc}}                       \\
             F(ab)F(c)   \ar@{=>}[rd]_{\varepsilon_{ab,c}}
          & &    F\big(a(bc)\big)        \\
          & F\big((ab)c\big) \ar@{=>}[ru]_{F(\phi)} &                  }\]
 Weak 2-categories and weak 2-functors between them form a
 category, which we denote by $\mathbf{W2Cat}$.

We can also talk about transformations between weak functors, and
modifications between them. The definitions are exactly the same
as the ones for the case of weak functors between strict
2-categories ($\S$\ref{SS:Transweak}).

Similarly, we have the categories $\mathbf{W2Gpd}$
 and $\mathbf{W2Gp}$. (We define  weak 2-functors between weak
 2-group(oid)s to be simply  weak 2-functors between the underlying weak
 2-categories.)
\subsection{Nerve of a weak 2-groupoid}{\label{SS:Weaknerve}}

To a weak 2-category $\mfC$ we can associate a nerve $N\mfC$.  The
simplices are exactly the same as in the case of strict 2-categories
(see Section~\ref{S:1}), except that now the commutativity of a
3-cell like the one on page \pageref{3cell} means
$[\phi_{f,g,m}][f\gamma][\beta]=[\alpha m][\delta]$. With some
careful, and not entirely easy, 2-diagram chasing, it can be shown
that taking nerves is a functor $N \: \mathbf{W2Gpd} \to
\mathbf{SSet2}$. This functor is indeed   an equivalence of
categories (Proposition~\ref{P:w2gpd}).

\vspace{0.1in}

\noindent{\em Some notation.} In what follows, the symbols $d_i$,
$s_j$  stand for the face and degeneracy  maps of a simplicial
set. The boundary of $\Delta^n$ is denoted by $\partial\Delta^n$.
The $k^{th}$ {\em horn } $\Lambda^n_k$ of $\Delta^n$ is the
  subcomplex of the simplicial set $\Delta^n$ generated by all but the $k^{th}$
   $(n-1)$-faces. If $x \in X_n$ is an n-simplex of a simplicial
   set $X$, we use the same notation for the corresponding map
   $x \: \Delta^n \to X$.

\vspace{0.1in}

\begin{lem}{\label{L:unique}}
  Let $X$ be a simplicial set in $\mathbf{SSet2}$. Let
  $x,y \: \Delta^3 \to X$ be 3-simplices that are equal on a
  horn $\Lambda^3_k$, for some $0 \leq k \leq 3$. Then $x=y$.
\end{lem}

\begin{proof}
   It is easy to see that $d_k(x)$ and $d_k(y)$ are equal
   on the boundary, and are homotopic relative to their equal
   boundary. So, by 2-minimality, $d_k(x)=d_k(y)$. It follows
   from Remark~\ref{R:minimal}.$\mathbf{b}$ that $x=y$.
\end{proof}

\begin{prop}{\label{P:w2gpd}}
   The functor $N \: \mathbf{W2Gpd} \to \mathbf{SSet2}$
   is an equivalence of categories.
\end{prop}

\begin{proof}[Sketch]
  We explain how to construct an inverse equivalence
  $\mathbf{SSet2} \to \mathbf{W2Gpd}$. Pick a simplicial set
  $X \in \mathbf{SSet2}$.  Define the weak 2-groupoid $\mfG_X$ as follows.
  The  set objects of $\mfG_X$ is $X_0$, and the set of arrows is $X_1$, with
  the source and target maps being $s=d_0$, $t=d_1$.
  Take the degenerate 1-simplices in $X_1$ for the identity arrows.
  To define composition of arrows, take
  two composable arrows $f$, $g \in \mfG_X$, and consider the induced
  map $\Lambda_2^1 \to X$.
    Since $X$ is Kan, we can extend this map to a
   map $I_{f,g} \: \Delta^2 \to X$; {\em choose} such an
   extension.
   Make sure to choose $I_{f,g}=s_0(f)$,
   respectively $I_{f,g}=s_1(g)$, when $g$, respectively $f$,
   is an identity arrow (that is, a degenerate 1-simplex in $X_1$).
   Define the composition $fg$ to be the restriction of $I_{f,g}$ to
   the $1^{st}$ face of $\Delta^2$. That is, $fg:=d_1(I_{f,g})$.

   For a pair of arrows $f$ and $g$ between objects $A$ and $B$,
   we define the set of 2-morphism from $f$ to $g$ to be the
   set of 2-simplices $\alpha \in X_2$ such that $d_2(\alpha)=f$, $d_1(\alpha)=g$,
   and $d_0(\alpha)=\id_B$. When $f=g$, we take the degenerate 2-simplex
   $s_1(f)$ for the identity 2-morphism from $f$ to itself .

   Using the fact that $X$ is Kan, plus the fact
   that $X$ satisfies the minimality assumption (for 2-simplices), one can
   show, for given objects $A$ and $B$, that there is a well-defined associative
   composition for
   2-morphisms, making $\hom(A,b)$ into a groupoid. The identity arrows in this hom-groupoid
   are the identity 2-morphisms defined in the previous
   paragraph.  One can also prove the following key fact:
    \begin{itemize}
      \item[$\blacktriangleright$]  Given 1-morphisms $f$, $g$ and $h$, the 2-morphisms
      from $fg$ to $h$ are in a natural bijection with the 2-simplices
      $\alpha \in X_2$ such that $f=d_2(\alpha)$, $g=d_0(\alpha)$ and
      $h=d_1(\alpha)$; when $h=fg$, the 2-isomorphism $\id \: fg \to h$
      corresponds to $I_{f,g}$ under this bijection.
    \end{itemize}

 Let us explain how to compose a 2-morphism with a 1-morphism.
 For example, suppose $f$ and $g$ are 1-morphisms from $A$ to
 $B$, $\alpha \: f \Rightarrow g$ is a 2-morphism  (as defined above),
 and $h$ is  a 1-morphism from $B$ to another object $C$.
 We define $\alpha h \: fh \Rightarrow gh$ as follows. Let $H \:
 \Lambda_2^3 \to X$ be the horn defined by $d_0(H)=s_0(h)$,
 $d_1(H)=I_{g,h}$, $d_3(H)=\alpha$. By the Kan property of
 $X$, together with Lemma~\ref{L:unique}, we obtain a unique
 extension of $H$ to a 3-simplex $\bar{H} \: \Delta^3 \to X$.
 The second face $d_2(\bar{H})$ of this 3-simplex determines,
 by~$\blacktriangleright$, a 2-morphism $fh \Rightarrow gh$,
 which we take to be $\alpha h$.

 Given three composable 1-morphisms $f$, $g$, and $h$, we define
 the associator $\phi_{f,g,h}$ as follows. Consider the horn
 $H \: \Lambda_1^3 \to X$ defined by $d_0(H)=I_{g,h}$,
 $d_2(H)=I_{f,gh}$, and $d_3(H)=I_{f,g}$. Let $\bar{H}$ be the
 (unique)
 extension of $H$ to a 3-simplex $\bar{H} \: \Delta^3 \to X$.
 We define $\phi_{f,g,h} \: (fg)h \Rightarrow f(gh)$ to be
 the 2-simplex
 $d_1(\bar{H})$; here, we are again using~$\blacktriangleright$.
 The following fact is the key to proving that $\mfG_X$ satisfies
 the axioms of a weak 2-groupoid, as well as to proving that the
 functor $X \mapsto \mfG_X$ is an inverse to the nerve functor.

  \begin{itemize}
      \item[$\blacktriangleleft\!\!\blacktriangleright$]
    Suppose we are given a hollow 3-simplex $\partial\Delta^3 \to X$
    whose faces are labeled as in the tetrahedron on page
    \pageref{nervediagram}. Then, this hollow 3-simplex
    can be extended to a (necessarily unique) 3-simplex
    $\Delta^3 \to X$
     if an only if $(f\gamma)(\beta)=(\alpha m)(\delta)$.
   \end{itemize}

Using $\blacktriangleleft\!\!\blacktriangleright$  we prove the
pentagon identity as follows. Let $\{a,b,c,d\}$ be a sequence of
composable maps. We will denote the sequence of vertices of these
arrows by $\{0,1,2,3,4\}$. Let $\Lambda \subset \Delta^4$ be the
2-skeleton of the standard 4-simplex. Consider the map $J \:
\Lambda \to X$ whose effect on 2-faces is given by:
\[\xymatrix@C=6pt@R=14pt@M=6pt{  & 1 \ar[dr]^b \ar@{}[d] | (0.7){I} &   \\
            0 \ar[rr]_{ab} \ar[ru]^a   &    &   2  } \ \
  \xymatrix@C=6pt@R=14pt@M=6pt{  & 1 \ar[dr]^{bc} \ar@{}[d] | (0.7){I} &   \\
            0 \ar[rr]_{a(bc)} \ar[ru]^a   &    &   3   } \ \
  \xymatrix@C=6pt@R=14pt@M=6pt{  & 1 \ar[dr]^{b(cd)} \ar@{}[d] | (0.7){I} &   \\
            0 \ar[rr]_{(a(b(cd))} \ar[ru]^a   &    &   4   } \ \
  \xymatrix@C=6pt@R=14pt@M=6pt{  & 2 \ar[dr]^c \ar@{}[d] | (0.7){\phi_{a,b,c}} &   \\
            0 \ar[rr]_{a(bc)} \ar[ru]^{ab}   &    &   3   }\]

\[\xymatrix@C=6pt@R=14pt@M=6pt{  & 2 \ar[dr]^{cd} \ar@{}[d] | (0.7){\phi_{a,b,cd}} &   \\
            0 \ar[rr]_{(a(b(cd))} \ar[ru]^{ab}   &    &   4   } \ \
  \xymatrix@C=6pt@R=14pt@M=6pt{  & 3 \ar[dr]^d \ar@{}[d] | (0.7){T} &   \\
            0 \ar[rr]_{(a(b(cd))} \ar[ru]^{a(bc)}   &    &   4} \ \
  \xymatrix@C=6pt@R=14pt@M=6pt{  & 2 \ar[dr]^c\ar@{}[d] | (0.7){I} &   \\
            1 \ar[rr]_{bc} \ar[ru]^b   &    &   3   }  \ \
  \xymatrix@C=6pt@R=14pt@M=6pt{  & 2 \ar[dr]^{cd} \ar@{}[d] | (0.7){I} &   \\
            1 \ar[rr]_{b(cd)} \ar[ru]^b   &    &   4   }\]

\[\xymatrix@C=6pt@R=14pt@M=6pt{  & 3 \ar[dr]^d \ar@{}[d] | (0.7){\phi_{b,c,d}} &   \\
            1 \ar[rr]_{b(cd)} \ar[ru]^{bc}   &    &   4 } \ \
  \xymatrix@C=6pt@R=14pt@M=6pt{  & 3 \ar[dr]^d \ar@{}[d] | (0.7){I} &   \\
            2 \ar[rr]_{cd} \ar[ru]^c   &    &   4  }\]

In the above diagrams, $I_{f,g}$ has been abbreviated to $I$
(which is meant to remind the reader of the identity 2-morphism),
and $T=(\phi_{a,bc,d})(a\phi_{b,c,d})$.

In fact, $J$ extends (uniquely) to the horn $H \: \Lambda_1^5 \to
X$. To see this, one just needs to verify that $J$ can be extended
to each of the faces $d_0$, $d_2$, $d_3$, and $d_4$, which is easily
done using~$\blacktriangleleft\!\!\blacktriangleright$. The horn $H$
can now be uniquely extended to a full 4-simplex $\bar{H} \:
\Delta^4 \to X$ using the Kan extension property and Remark
\ref{R:minimal}.$\mathbf{b}$. This in particular implies that the, a
priori hollow, 3-cell $d_1(\bar{H})$ has a filling.
Using~$\blacktriangleleft\!\!\blacktriangleright$, this translates to
 \[(\phi_{a,b,c}d)(\phi_{a,bc,d})(a\phi_{b,c,d})
                          =(\phi_{ab,c,d})(\phi_{a,b,cd}),\]
which is the desired pentagon identity. This completes the (sketch
of the) proof that $\mfG_X$ is a weak 2-groupoid.

Let us now verify the functoriality of our construction. Namely,
given a simplicial map $F \: X \to Y$ in $\mathbf{SSet2}$, we
describe the induced weak functor $\mfG_F \: \mfG_X \to \mfG_Y$. The
effect of $\mfG_F$ on objects, 1-morphisms and 2-morphisms is
defined in the obvious way. Given a pair $f$ and $g$ of composable
morphisms in $\mfG_X$, we define $\epsilon_{f,g}$ to be the
2-morphism in $\mfG_Y$ associated to the 2-simplex $F(I_{f,g})$; we
are using~$\blacktriangleright$ again.

The fact that the above two functors are inverse to each other is
not difficult to check using~$\blacktriangleright$
and~$\blacktriangleleft\!\!\blacktriangleright$.
\end{proof}

Definition~\ref{D:homotopygroup} and Lemma~\ref{L:weakweak} remain
valid for weak maps between weak 2-groupoids.

\subsection{Strictifying weak 2-groupoids}{\label{SS:Strictify}}

Most of the 2-groups that appear in nature are weak. There are,
however, general methods to strictify  weak 2-groups.  We briefly
address this problem after the following example. We point out that
the question of strictifying weak 2-groups, or more generally
bicategories, has been studied by many authors (see for instance,
\cite{M-P}, \cite{Laplaza}, \cite{Baez}, \cite{Lack}), but our point
of view is slightly different, in that our main emphasis is on
showing that the strictifications that we consider do respect
derived mapping spaces. In other words, our main interest is in
using strictifications (or weakenings, for that matter) as tools for
computing derived mapping spaces.

\begin{ex}{\label{E:cocycles}}
 Let $A$ be an abelian group, and $\Gamma$ an arbitrary group.
 Consider the groupoid $\mathfrak{E}$ of central extensions of
 $\Gamma$ by $A$.
 There is a natural multiplication on
 $\mathfrak{E}$ which makes it into a weak 2-group.
 It is defined as follows. Take two objects
  \[\epsilon_1: 1 \arr{} A \arr{\alpha} E_1 \arr{} \Gamma \arr{} 1\]
  \[ \epsilon_2: 1 \arr{} A \arr{\beta} E_2 \arr{} \Gamma \arr{} 1\]
in $\mathcal{E}(\Gamma,A)$.  The sum
 $\epsilon_1+\epsilon_2$ is the sequence
    \[1 \arr{} A \arr{(\alpha,\beta^{-1})} E_1\oux{\Gamma}{A}E_2 \arr{} \Gamma \arr{} 1,\]
where $\beta^{-1}$ stands for the point-wise inversion of $\beta$.
The inverse of $\epsilon_1$ is given by
    \[ \epsilon_1^{-1}: 1 \arr{} A \arr{\alpha^{-1}} E_1 \arr{} \Gamma \arr{} 1.\]

 The set of isomorphism classes
 of elements in $\mathcal{E}(\Gamma,A)$ is in natural bijection
 with $H^2(\Gamma,A)$, and the automorphism group of an object in
$\mathfrak{E}$  is naturally
  isomorphic to $H^1(\Gamma,A)$.
 Therefore, we have $\pi_1\mathfrak{E}=H^2(\Gamma,A)$, and
 $\pi_2\mathfrak{E}=H^1(\Gamma,A)$.

 We can give a strict model for the weak 2-group $\mathfrak{E}$
 using cocycles  as follows. Let $C^2(\Gamma,A)$
 be the group of 2-cocycles from $\Gamma$ to $A$.
 Let $A^{\Gamma}$ be the group of (set) maps from $\Gamma$ to $A$.
 Let $C^2(\Gamma,A)$ act trivially
 on $A^{\Gamma}$. Form the crossed module
 $[\partial \: A^{\Gamma} \to C^2(\Gamma,A)]$, where
 $\partial$ is the boundary map.  This crossed module is a strict model
 for $\mathfrak{E}$.
  In fact, this example can be made to work in the case where $A$
  is an arbitrary $\Gamma$-module, and $\mathfrak{E}$  is the
  group of such extensions for which the induced action of
  $\Gamma$  on $A$ coincides with the given one.
\end{ex}

Now we consider the general strictification problem for weak
2-groupoids and weak maps. We have a hierarchy of 2-groupoids as in
the following commutative diagram:
   \[\xymatrix@=12pt@M=10pt{ \mathbf{2Gpd}
      \ar@{}|-{\subset}[r]\ar[rd]_(0.39){N}  &  \mathbf{2Gpd}_w \ar@{}|-{\subset}[r]
                                \ar[d]_{N}& \mathbf{W2Gpd}  \ar[dl]^(0.43){N}  \\
              &  \mathbf{SSet}  &       }\]
All the maps in this diagram respect (and reflect) weak
equivalences. So we can pass to  homotopy categories.

\begin{prop}{\label{P:equivalent1}}
  We have a commutative diagram of equivalences of categories
      \[\xymatrix@=12pt@M=8pt{ \Ho(\mathbf{2Gpd})
   \ar@{}|-{\subset}[r]_(0.49){\sim}\ar[rd]_(0.40){\sim}  &
                 \Ho(\mathbf{2Gpd}_w) \ar@{}|-{\subset}[r]_(0.49){\sim}
                   \ar[d]_{\sim}& \Ho(\mathbf{W2Gpd})  \ar[dl]^(0.44){\sim}  \\
              &  \mathbf{2Types}  &       }\]
 Here $\mathbf{2Types}$ stands for the full subcategory of $\Ho(\mathbf{SSet})$
consisting of all simplicial sets $X$ with\/ $\pi_iX=0$, $i\geq 3$.
In all three cases, the downward arrows are given by $N$, and their
inverse equivalences are given by $W$.
\end{prop}

\begin{proof} We have to prove that $W$ provides an inverse
equivalence to all of the downward arrows. For the leftmost arrow,
this is part of Theorem~\ref{T:QuillenEq}. Let us   prove the
statement for the rightmost arrow. We have to show that $W \:
\mathbf{2Types}  \to \Ho(\mathbf{W2Gpd})$ is an inverse equivalence
to  $N \: \Ho(\mathbf{W2Gpd}) \to \mathbf{2Types}$. The fact that
$N\circ W \: \mathbf{2Types} \to \mathbf{2Types}$ is equivalent to
$\id_{\mathbf{2Types}}$ is already part of Theorem
\ref{T:QuillenEq}, because $W \: \mathbf{2Types}  \to
\Ho(\mathbf{W2Gpd})$ factors through $\Ho(\mathbf{2Gpd})$. We now
show that for every weak 2-groupoid $\mfG$ there is a weak
equivalence $\mfG \to WN(\mfG)$, natural in $\mfG$. Since $N \:
\mathbf{W2Gpd} \to \mathbf{SSet}$ is fully faithful and reflects
weak equivalences, it is enough to construct a weak equivalence
 $N\mfG \to NWN(\mfG)$. But this is easy because
 we know that for every simplicial set
 $X$ with $\pi_iX=0$ for $i\geq 3$, for example $X=N\mfG$,
 there is a natural map $X \to NW(X)$, namely the unit
 of the adjunction of Theorem~\ref{T:QuillenEq}, and that this map
 is a weak equivalence (Theorem~\ref{T:QuillenEq}.$\mathbf{iii}$).

The proof for the case of $\Ho(\mathbf{2Gpd}_w)$ is similar.
\end{proof}

The same discussion applies to the case where everything is pointed.
In particular, we have a diagram
  \[\xymatrix@=12pt@M=10pt{ \mathbf{2Gp}
   \ar@{}|-{\subset}[r]\ar[rd]_(0.39){N}  &  \mathbf{2Gp}_w \ar@{}|-{\subset}[r]
                              \ar[d]_{N}& \mathbf{W2Gp}  \ar[dl]^(0.43){N}  \\
              &  \mathbf{SSet}_*  &       }\]
which passes to the homotopy category, and we obtain the following
proposition whose proof is similar to that of Proposition
\ref{P:equivalent1}.

\begin{prop}{\label{P:equivalent2}}
  We have a commutative diagram of equivalences of categories
     \[\xymatrix@=12pt@M=8pt{ \Ho(\mathbf{2Gp})
   \ar@{}|-{\subset}[r]_(0.49){\sim}\ar[rd]_(0.40){\sim}  &
                 \Ho(\mathbf{2Gp}_w) \ar@{}|-{\subset}[r]_(0.49){\sim}
                   \ar[d]_{\sim}& \Ho(\mathbf{W2Gp})  \ar[dl]^(0.44){\sim}  \\
              &  \mathbf{Con2Types}_*  &       }\]
 Here $\mathbf{Con2Types}_*$ is the full subcategory of $\Ho(\mathbf{SSet}_*)$
consisting of pointed connected simplicial sets $X$ with\/
$\pi_iX=0$, $i\geq 3$. In all three cases, the downward arrows are
given by $N$, and their  inverse equivalences are given by $W$.
\end{prop}

Propositions~\ref{P:equivalent1} and~\ref{P:equivalent2} can in fact
be strengthened to statements about simplicial mapping spaces.
This means that all these functors
 respect derived  mapping spaces, up to natural
homotopy equivalences. Let us explain what we mean by ``derived.''
If we are in $\mathbf{2Gpd}$, $\mathbf{2Gpd}_*$ or $\mathbf{2Gp}$,
derived simply means  working with $\Mapbf$ rather than $\homs$. If
we are in $\mathbf{SSet}$, to compute derived mapping spaces we need
to first make a fibrant replacement on the target; observe, however,
that the images of the nerve functors always land in the fibrant
part, so in such cases derived mapping spaces are the same as the
usual mapping spaces $\Hombf$. In $\mathbf{2Gpd}_w$,
$\mathbf{W2Gpd}$, and in their pointed versions,  derived mapping
space simply means the  hom-2-groupoid $\homw$.

For instance, if we are given weak 2-groups $\mfH$ and $\mfG$,
then the hom-2-groupoid $\homw(\mfH,\mfG)$ is naturally equivalent
to the derived hom-2-groupoid $\Mapbf(\mfH^s,\mfG^s)$, where
$\mfG^s:=WN\mfG$ is the ``strictification'' of $\mfG$. This
hom-2-groupoid is in turn naturally equivalent to
$\homw(\mfH^s,\mfG^s)$, where the latter is computed in
$\mathbf{2Gpd}_w$.

\begin{prop}{\label{P:weakversion}}
 Suppose $\mfG$ and $\mfH$ are weak 2-groupoids.
 \begin{itemize}
  \item[$\mathbf{i.}$] Let $\mfH' \to \mfH$ and $\mfG \to \mfG'$
    be equivalences of weak 2-groupoids. Then the induced
     functor $\homw(\mfH,\mfG) \to \homw(\mfH',\mfG')$ is an
     equivalence of weak 2-groupoids.

  \item[$\mathbf{ii.}$] We have a natural equivalence
    of weak 2-groupoids
     \[\Mapbf(\mfH',\mfG')\simeq\homw(\mfH,\mfG),\]
    where  $\mfG \risom \mfG'$  and $\mfH' \risom \mfH$
    are arbitrary
    strictifications of $\mfH$ and $\mfG$; no requirement on
    these maps being cofibration or fibration.

  \item[$\mathbf{iii.}$] We have natural equivalences of
  simplicial sets
   \[\Hombf(N\mfH,N\mfG)\simeq N\Mapbf(\mfH',\mfG')
                             \simeq N\homw(\mfH,\mfG),\]
     where $\mfG'$ and $\mfH'$ are as in part ($\mathbf{ii}$).
 \end{itemize}
\end{prop}

\begin{proof}
  Part ($\mathbf{i}$) is proved by an argument similar to the one
  used in the proof of Proposition~\ref{P:mappingspace2}.  Part ($\mathbf{ii}$)
  follows from  part ($\mathbf{i}$) and
  Proposition~\ref{P:mappingspace2}. The second equivalence of  ($\mathbf{iii}$)
  follows from ($\mathbf{ii}$). For the first equivalence, note
  that,   since the nerve of a 2-groupoid is a fibrant simplicial set,
  we have $\Hombf(N\mfH,N\mfG)\simeq \Hombf(N\mfH',N\mfG')$.
  Now use Corollary~\ref{C:mappingspace}.
\end{proof}

\begin{cor}{\label{C:lift2}}
  Suppose $\mfG$ and $\mfH$ are  weak 2-groupoids. Then, for every
  map $f \in [\mfH,\mfG]_{\mathbf{W2Gpd}}$,  there is a weak map
  $\tilde{f} \: \mfH \to \mfG$,
  unique up to transformation, which induces $f$ in the homotopy
  category.
\end{cor}

\begin{proof}
  Using Proposition~\ref{P:weakversion}, the proof of
  Corollary~\ref{C:lift} works in this case as well.
\end{proof}

\begin{prop}{\label{P:pointed}}
  The pointed versions of Proposition~\ref{P:weakversion}
  and Corollary~\ref{C:lift2} are true. In particular, we have similar
  statements for 2-groups and crossed modules.
\end{prop}

\providecommand{\bysame}{\leavevmode\hbox
to3em{\hrulefill}\thinspace}
\providecommand{\MR}{\relax\ifhmode\unskip\space\fi MR }
\providecommand{\MRhref}[2]{%
  \href{http://www.ams.org/mathscinet-getitem?mr=#1}{#2}
} \providecommand{\href}[2]{#2}


\begin{thebibliography}{10}


\bibitem[BaLa]{Baez} J.~Baez and A.~Lauda, Higher
dimensional algebra {V}: 2-groups,  \emph{Theory Appl. Categ.}
\textbf{12}  (2004), 423--491.

\bibitem[Ban]{Bangor} Webpage of the Bangor Reserach Group:\\
\textsf{http://www.informatics.bangor.ac.uk/public/mathematics/research/cathom/cathom2.html}

\bibitem[Ba]{Baues} H-J.~Baues, \emph{Combinatorial homotopy and
$4$-dimensional complexes}, de Gruyter Expositions in Mathematics \textbf{2},
Walter de Gruyter, Berlin, 1991.

\bibitem[Be]{B} J.~B{\'e}nabou, Introduction to bicategories,
\emph{Reports of the Midwest Category Seminar},
Lecture Notes in Mathematics, vol.~47, 1--77,
Springer-Verlag, New York, 1967.

\bibitem[Br]{Brown} K.~Brown, \emph{Cohomology of Groups},
Graduate Texts in Mathematics 87,
Springer-Verlag, New York, 1994.


\bibitem[BrGo]{BrGo} R.~Brown, M.~Golasinski, \emph{A model structure for the homotopy theory of crossed complexes},
Cahiers Topologie G\'eom. Diff\'erentielle Cat\'eg.  \textbf{30}
(1989),  no. 1, 61--82.

\bibitem[Gr]{Gray} J.W.~Gray, \emph{Formal category theory:
adjointness for $2$-categories},
  Lecture Notes in Mathematics, vol.~391, Springer-Verlag, Berlin-New York,
  1974.


\bibitem[GoJa]{Jardine} P.~G.~Goerss and J.~F.~Jardine,
\emph{Simplicial homotopy theory}, Progress in
  Mathematics 174, Birkh\"{a}user Verlag, Basel, 1999.

\bibitem[La02]{Lack2} S.~Lack, A Quillen model structure for 2-categories,
   \emph{$K$-Theory} \textbf{26} (2002), 171--205.

\bibitem[La04]{Lack} S.~Lack, A Quillen model structure for
bicategories, \emph{$K$-Theory} \textbf{33}  (2004),  185--197.

\bibitem[Lap]{Laplaza} M.~L.~Laplaza, Coherence for
categories with group structure: an alternative approach,
\emph{J. Algebra} \textbf{84} (1983),  no. 2, 305--323.

\bibitem[Lo]{Loday} J-L.~Loday, Spaces with finitely many
nontrivial homotopy groups,
\emph{J. Pure Appl. Algebra} \textbf{24} (1982), no.~2, 179--202.

\bibitem[McPa]{M-P} S.~MacLane, R.~Par\'e, Coherence for bicategories
 and indexed categories,  \emph{J. Pure Appl. Algebra}  \textbf{37}  (1985),
   no.~1, 59--80.

\bibitem[McWh]{M-W} S.~MacLane, J.~H.~C.~Whitehead, On the
$3$-type of a complex,
  \emph{Proc. Nat. Acad. Sci. U.S.A.} \textbf{36} (1950), 41--48.

\bibitem[Ma]{May} J.~P.~May, \emph{Simplicial objects in algebraic
topology}, Chicago Lectures in
  Mathematics, University of Chicago Press, Chicago, 1992.


\bibitem[MoSe]{M-S} I.~Moerdijk and J.~Svensson, Algebraic
classification of equivariant
  homotopy 2-types, \emph{J. Pure Appl. Algebra} \textbf{89} (1993), 187--216.

\bibitem[No]{Maps} B.~Noohi, On weak maps between 2-groups,
  preprint, arXiv:math/0506313v2 [math.CT].

\bibitem[To]{Tonk} A.~Tonks, On the Eilenberg-Zilber theorem for crossed complexes,
\emph{J. Pure Appl. Algebra}  \textbf{179}  (2003), 199--220.

\bibitem[Wh]{W} J.~H.~C.~Whitehead, Combinatorial homotopy
{I}, \emph{Bull. Amer. Math. Soc.}
  \textbf{55} (1949), 213--245.

\end{thebibliography}
\end{document}